# Optimal Impact Strategies for Asteroid Deflection

Massimiliano Vasile[*] and Camilla Colombo[†]

*Glasgow University, Glasgow, UK, G128QQ*

This paper presents an analysis of optimal impact strategies to deflect potentially dangerous asteroids. In order to compute the increase in the minimum orbit intersection distance of the asteroid due to an impact with a spacecraft, simple analytical formulae are derived from proximal motion equations. The proposed analytical formulation allows for an analysis of the optimal direction of the deviating impulse transferred to the asteroid. This ideal optimal direction cannot be achieved for every asteroid at any time, therefore an analysis of the optimal launch opportunities for deviating a number of selected asteroids was performed through the use of a global optimization procedure. The results in this paper demonstrate how the proximal motion formulation has very good accuracy in predicting the actual deviation and can be used with any deviation method since it has general validity. Furthermore the characterization of optimal launch opportunities shows how a significant deviation can be obtained even with a small spacecraft.

## Nomenclature

$\mathbf{A}_{MOID}$ = matrix form of the proximal motion equations

$a$ = semi-major axis of the nominal orbit, km

$b$ = semi-minor axis of the nominal orbit, km

$b^*$ = impact parameter, km

$e$ = eccentricity of the nominal orbit

---

[*] Lecturer Ph.D., Department of Aerospace Engineering, James Watt (South Building), m.vasile@eng.gla.ac.uk, AIAA Member.
[†] Ph.D. Candidate, Department of Aerospace Engineering, James Watt (South Building), camilla.colombo@strath.ac.uk, AIAA Student Member.

| | | |
|---|---|---|
| $e_r$ | = | relative error |
| $g_0$ | = | standard free-fall (9.81 m/s²) |
| $\mathbf{G_d}$ | = | matrix form of the Gauss' equations |
| $h$ | = | angular momentum of the nominal orbit, km²/s |
| $\hat{\mathbf{h}}$ | = | unit vector in the direction of angular momentum |
| $i$ | = | inclination of the nominal orbit, deg |
| $I_{sp}$ | = | specific impulse of the spacecraft engine, s |
| $J$ | = | objective function of the optimization problem |
| $M$ | = | mean anomaly of the nominal orbit, deg |
| $m$ | = | mass, kg |
| $n$ | = | angular velocity of the nominal orbit, s⁻¹ |
| $\hat{\mathbf{n}}$ | = | unit vector in the normal direction in the orbital plane |
| $p$ | = | semilatum rectum of the nominal orbit, km |
| $r$ | = | norm of the position on the nominal orbit, km |
| $T$ | = | transition matrix |
| $T_{NEO}$ | = | period of the asteroid, d |
| $\hat{\mathbf{t}}$ | = | unit vector in the tangential direction in the orbital plane |
| $t_d$ | = | time at which the impulsive deviation is performed, s or d |
| $t_{MOID}$ | = | time at the Minimum Orbit Interception Distance point, s or d |
| $t_w$ | = | warning time, d |
| $v$ | = | orbital velocity, km/s |
| $\delta\boldsymbol{\alpha}$ | = | orbital elements difference |
| $\Delta\mathbf{r}$ | = | vector distance of the asteroid from Earth at the Minimum Orbit Interception Distance, km |
| $\Delta r$ | = | norm of $\Delta\mathbf{r}$, km |
| $\delta\mathbf{r}$ | = | deviation vector in the Hill coordinate frame, km |
| $\delta r$ | = | norm of $\delta\mathbf{r}$, km |
| $\Delta\mathbf{v}$ | = | relative velocity vector, km/s |
| $\delta\mathbf{v}$ | = | impulsive maneuver vector, km/s |
| $\delta\hat{\mathbf{v}}_{opt}$ | = | optimal impulse unit vector |
| $\mu$ | = | gravitational constant, km$^3$/s$^2$ |
| $\mathbf{v}$ | = | eigenvectors of the transition matrix |
| $\Omega$ | = | argument of the ascending node of the nominal orbit, deg |
| $\omega$ | = | argument of the perigee of the nominal orbit, deg |
| $\theta$ | = | true anomaly on the nominal orbit, deg |
| $\theta^*$ | = | argument of latitude on the nominal orbit, deg |

*Subscripts*

| | | |
|---|---|---|
| $(\ )_d$ | = | value calculated at the instant the deviation maneuver is performed |
| $(\ )_{MOID}$ | = | value calculated at the Minimum Orbit Interception Distance |

## I.   Introduction

The Earth, as most of the planets of the Solar System, from its formation up to recent times, has experienced a strong interaction with minor celestial bodies such as asteroids and comets, many times resulting in dramatic impact events which have had serious effects on the geological, climate and global evolution of our planet. Though it is not possible to assess how many times the Earth during its history has been hit by one of these objects, the some 170 impact craters recognized around the globe give a scientific evidence that such astronomical events have repeatedly happened in the past, sometimes with global and dramatic consequences on the planet, as in the case of the Chicxulub crater in the Yucatan, which is estimated to have been caused by a celestial body of about ten km in size, and which is today widely acknowledged as the reason for the dinosaurs' extinction some 65 million years ago. In fact the potential impact of such a large and massive object, though statistically unlikely, having a probability of one event over millions of years, would certainly pose a critical threat to most of the population of the planet. It is estimated that there are more than 1000 objects having a size larger than 1 km. On the other hand if these large size bodies are extremely rare, it should be noticed that objects of 40 m in diameter, which is considered the critical threshold above which the Earth's atmosphere is no longer disintegrating an object, are estimated to be more than one million in number with a statistical frequency of impact of one hundred years or even less. An example of such an event happened at the beginning of the twentieth century in Siberia where an object of few tens of meters, though disintegrating before hitting the ground, produced a large devastation of many square kilometers of the Siberian forest. It is quite clear what the effect of such an impact would be if it happened in a densely populated area.

The recent discovery of Apophis, for which an impact in 2029 has been definitively excluded, through a passage in a keyhole during the 2029 fly-by could still lead to a close approach in 2036, has widely drawn the attention of public and media on the issue of potentially hazardous objects, and consequently of the technological and detection capabilities that nations have in order to implement a mitigation and prevention policy. Many research centers and groups exist around the world currently using ground based observation methodologies in order to keep track of risky objects together with a continuous update of the risk assessment related to each potential hazardous object. At the same time space agencies have started widening their scope to comets and asteroids not only for scientific reasons, but also in the frame of developing technological capabilities required in the case of an object should pose a serious threat to the Earth.

The European Space Agency in particular is now assessing the feasibility of the Don Quijote mission[1], due to launch in the first half of next decade, which is intended to impact a spacecraft with a high relative velocity onto an asteroid and measure its deflection. Should this mission fly, this would be the first technological demonstration of our capability to deviate an asteroid if needed.

Prevention strategies against a potential hazardous object in collision route with the Earth usually consider a change in momentum of the asteroid, with a consequent variation in the semi-major axis which results in an increase of the Minimum Orbit Intersection Distance (MOID), between the Earth and the object. Several different strategies have been considered to achieve this goal; among them the simplest one is the kinetic impact. In fact, as will be shown in this paper, effective kinetic impacts resulting in a variation of the MOID even of thousand of kilometers seem to be already achievable with the current launch technology with a relatively small spacecraft, provided that the time difference between the momentum change and the potential Earth impact is large enough.

In this paper we focus our attention on the analysis of optimal impact strategies for the deflection of asteroids. A simple analytical expression is derived to compute the displacement of the position of an asteroid at the Minimum

Orbit Interception Distance point, after the impact with a small spacecraft. This analytical formulation makes use of proximal motion equations expressed as a function of orbital elements[2]. These equations, already available in the literature, provide, with very little modifications, a very simple and general means to compute with a good accuracy the variation of the MOID without the need for further analytical developments.

The approach presented in this paper is an extension and a generalization of the methodologies proposed in previous works[3-5] where formulae were devised to compute the deviation due to a variation of the orbital mean motion, i.e. due to an action applied along the direction of the motion of the asteroid. A more general approach was already proposed by Conway[6] who used Lagrange coefficients expressed through universal formulae to analytically propagate forward in time the result of an impulsive maneuver. By that, he derived an elegant way to determine the near-optimal direction in which the impulse should be given. A similar approach was taken by Park and Ross[7], who performed a constrained optimization to compute the optimal impulsive deflection maneuver; in a subsequent work[8] they added the gravitational effects of Earth obtaining a more accurate estimation of the optimal impulse.

In this paper, near-optimal directions for deviation impulses are derived using a simple restricted two-body dynamic model. The gravitational effect of the Earth is accounted for by looking at the obtained deviation on the b-plane. The accuracy of the result is then assessed using a numerical propagation of the post deviation conditions with a full three-body dynamic model including the Sun and the Earth.

Since ideal optimal deflection conditions cannot always be achieved, a characterization of optimal mission opportunities is performed for a restricted group of selected asteroids over a very wide range of possible launch dates. It was decided to look only for mission options with a relatively low transfer time therefore only direct transfers and transfers with one single swing-by of Venus were considered. Longer sequences of swing-bys, though improving the deviation, would imply a longer term planning and more complex operations. The search for different transfer options was performed with a particular global optimization procedure based on an automatic branch and prune of the solution space combined with an agent-based search technique. This methodology has been proven to be an effective way of looking for families of optimal solutions in the case of complex trajectory design problems[9,10]. The result of this analysis demonstrates how, with a small spacecraft and very simple transfer strategies, it is possible to obtain considerable deviations for most of the threatening asteroids.

## II. Maximum Deviation Problem

Given the Minimum Orbit Interception Distance from the Earth for a generic Near Earth Object (NEO), the objective is to maximize the MOID by applying an impulsive action at a certain time $t_d$. The impulse acts as a perturbation on the orbit of the NEO and its new orbit can be considered proximal to the unperturbed one. If $\theta_{MOID}$ is the true anomaly of the NEO at the MOID along the unperturbed orbit, and $\theta^*_{MOID} = \theta_{MOID} + \omega$ the corresponding argument of latitude, we can write the variation of the position of the NEO, after deviation, with respect to its unperturbed position by using the proximal motion equations[2]:

$$\delta s_r \approx \frac{r}{a}\delta a + \frac{ae\sin\theta_{MOID}}{\eta}\delta M - a\cos\theta_{MOID}\delta e$$

$$\delta s_\theta \approx \frac{r}{\eta^3}(1+e\cos\theta_{MOID})^2 \delta M + r\delta\omega + \frac{r\sin\theta_{MOID}}{\eta^2}(2+e\cos\theta_{MOID})\delta e + r\cos i\delta\Omega \quad (1)$$

$$\delta s_h \approx r\left(\sin\theta^*_{MOID}\delta i - \cos\theta^*_{MOID}\sin i\delta\Omega\right)$$

where $\delta s_r$, $\delta s_\theta$ and $\delta s_h$ are the displacements in the radial, transversal and perpendicular to the orbit plane directions respectively, such that $\delta \mathbf{r} = [\delta s_r \quad \delta s_\vartheta \quad \delta s_h]^T$, and $\eta = \sqrt{1-e^2}$. The variation of the orbital parameters $a$, $e$, $i$, $\Omega$ and $\omega$ are computed through Gauss' planetary equations[11] considering an instantaneous change in the NEO velocity vector $\delta \mathbf{v} = [\delta v_t \quad \delta v_n \quad \delta v_h]^T$:

$$\delta a = \frac{2a^2 v}{\mu}\delta v_t$$

$$\delta e = \frac{1}{v}\left[2(e+\cos\theta_d)\delta v_t - \frac{r}{a}\sin\theta_d \delta v_n\right]$$

$$\delta i = \frac{r\cos\theta^*_d}{h}\delta v_h$$

$$\delta\Omega = \frac{r\sin\theta^*_d}{h\sin i}\delta v_h \quad (2)$$

$$\delta\omega = \frac{1}{ev}\left[2\sin\theta_d \delta v_t + \left(2e+\frac{r}{a}\cos\theta_d\right)\delta v_n\right] - \frac{r\sin\theta^*_d \cos i}{h\sin i}\delta v_h$$

$$\delta M_{t_d} = -\frac{b}{eav}\left[2\left(1+\frac{e^2 r}{p}\right)\sin\theta_d \delta v_t + \frac{r}{a}\cos\theta_d \delta v_n\right]$$

the above variation on $M$ takes into account only the instantaneous change of the orbit geometry at time $t_d$. On the other hand due to the change in the semi-major axis we have a variation of the mean motion $n$ and therefore a change in the mean anomaly at the time of the MOID given by:

$$\delta M_n = \delta n \left( t_{MOID} - t_d \right) = \delta n \Delta t \qquad (3)$$

where $t_{MOID}$ is the time at the MOID along the orbit of the NEO and $\delta n = \sqrt{\dfrac{\mu}{(a+\delta a)^3}} - \sqrt{\dfrac{\mu}{a^3}}$. $\Delta t$ is the time-to-impact, defined as $t_{MOID} - t_d$. Eq. (3) takes into account the phase shifting between the Earth and the NEO. The total variation in the mean anomaly between the unperturbed and the proximal orbit is therefore:

$$\delta M = -\frac{b}{eav}\left[ 2\left(1 + \frac{e^2 r}{p}\right)\sin\theta_d \delta v_t + \frac{r}{a}\cos\theta_d \delta v_n \right] + \delta n \Delta t \qquad (4)$$

Now if $\Delta \mathbf{r} = \begin{bmatrix} \Delta s_r & \Delta s_\theta & \Delta s_h \end{bmatrix}^T$ is the vector distance of the asteroid from the Earth at the MOID, and $\delta \mathbf{r} = \begin{bmatrix} \delta s_r & \delta s_\theta & \delta s_h \end{bmatrix}^T$ is the variation given by Eqs. (1) at $t_{MOID}$, then the objective function for the maximum deviation problem is the following:

$$J = \left(\Delta s_r + \delta s_r\right)^2 + \left(\Delta s_\theta + \delta s_\theta\right)^2 + \left(\Delta s_h + \delta s_h\right)^2 \qquad (5)$$

The proposed formulation provides a very fast analytical way of computing the variation of any asteroid's orbit due to any impulsive deviation action. It is both an extension and a generalization of other approaches[3,4,5] that were instead limited to actions in the direction of the unperturbed velocity of the asteroid. Compared to more general methods that are analytically propagating the perturbed trajectory by using the Lagrange coefficients[7], the proposed approach does not require any solution of the time equation for every variation of the orbit of the asteroid and therefore it is less computationally expensive. On the other hand, it is conceptually and computationally equivalent to those approaches[6] that are analytically propagating only the variation of the position and velocity of the asteroid by using the fundamental perturbation matrix[11], although the accuracy of the two formulations was not compared. Conversely, the benefit of using proximal motion equations expressed in orbital elements is the explicit relationship between the components of the perturbing action and the variation of the geometric characteristics of the orbit of the asteroid.

### A. Maximum Deviation Strategies

By combining Eqs. (1) and Eqs.(2), it is possible to compute for each time $t_d$ the transition matrix that links $\delta \mathbf{v}$ at $t_d$ to $\delta \mathbf{r}$ at $t_{MOID}$. In order to make explicit the dependence on the impulse components in each of the equations, Eq. (3) has to be rewritten as a function of $\delta a$ as follows:

$$\delta M_n = \delta n \Delta t = -\frac{3}{2}\frac{\sqrt{\mu}}{a^{\frac{5}{2}}}\Delta t \delta a \qquad (6)$$

If now Eq. (4) is incorporated into system (2) and along with Eqs. (1), are written in matrix form, we have:

$$\begin{cases} \delta \mathbf{r}(t_{MOID}) = \mathbf{A}_{MOID}\delta \mathbf{a}(t_d) \\ \delta \mathbf{a}(t_d) = \mathbf{G}_d \delta \mathbf{v}(t_d) \end{cases} \Rightarrow \delta \mathbf{r}(t_{MOID}) = \mathbf{A}_{MOID}\mathbf{G}_d \delta \mathbf{v}(t_d) = \mathbf{T}\delta \mathbf{v}(t_d) \qquad (7)$$

where $\delta \mathbf{a}(t_d) = [\delta a \;\; \delta e \;\; \delta i \;\; \delta \Omega \;\; \delta \omega \;\; \delta M]^T$ is the orbital element difference, $\mathbf{A}_{MOID}$ is the matrix:

$$\mathbf{A}_{MOID}^T = \begin{bmatrix} \frac{r_{MOID}}{a} - \frac{3}{2}\frac{e\sin\theta_{MOID}}{\eta}\frac{\sqrt{\mu}}{a^{\frac{3}{2}}}\Delta t & -\frac{3}{2}\frac{r_{MOID}}{\eta^3}(1+e\cos\theta_{MOID})^2 \frac{\sqrt{\mu}}{a^{\frac{5}{2}}}\Delta t & 0 \\ -a\cos\theta_{MOID} & \frac{r_{MOID}\sin\theta_{MOID}}{\eta^2}(2+e\cos\theta_{MOID}) & 0 \\ 0 & 0 & r_{MOID}\sin\theta^*_{MOID} \\ 0 & r_{MOID}\cos i & -r_{MOID}\cos\theta^*_{MOID}\sin i \\ 0 & r_{MOID} & 0 \\ \frac{ae\sin\theta_{MOID}}{\eta} & \frac{r_{MOID}}{\eta^3}(1+e\cos\theta_{MOID})^2 & 0 \end{bmatrix} \qquad (8)$$

and $\mathbf{G}_d$ is the matrix:

$$\mathbf{G}_d = \begin{bmatrix} \frac{2a^2 v_d}{\mu} & 0 & 0 \\ \frac{2(e+\cos\theta_d)}{v_d} & -\frac{r_d}{av_d}\sin\theta_d & 0 \\ 0 & 0 & \frac{r_d \cos\theta^*_d}{h} \\ 0 & 0 & \frac{r_d \sin\theta^*_d}{h\sin i} \\ \frac{2\sin\theta_d}{ev_d} & \frac{2e+\frac{r_d}{a}\cos\theta_d}{ev_d} & -\frac{r_d \sin\theta^*_d \cos i}{h\sin i} \\ -\frac{b}{eav_d}2\left(1+\frac{e^2 r_d}{p}\right)\sin\theta_d & -\frac{b}{eav_d}\frac{r_d}{a}\cos\theta_d & 0 \end{bmatrix} \qquad (9)$$

The subscript indices, $MOID$ and $d$, indicate that the matrices are calculated respectively at $t_{MOID}$ and $t_d$. As suggested by Conway[6], in order to maximize $\|\delta \mathbf{r}(t_{MOID})\| = \max(\mathbf{T}\delta \mathbf{v}(t_d))$, the associated quadratic form $\mathbf{T}^T\mathbf{T}$ has to be maximized by choosing an impulse vector $\delta \mathbf{v}(t_d)$ parallel to the eigenvector $\mathbf{v}$ of $\mathbf{T}^T\mathbf{T}$, conjugated to the maximum eigenvalue. Fig. 1a and b represent the components of the optimal impulse unit vector $\delta \hat{\mathbf{v}}_{opt}(t_d)$, projected onto the $\{t,n,h\}$ reference frame (where $\hat{\mathbf{t}}$ is the along the direction of motion, $\hat{\mathbf{h}}$ is the direction of the angular momentum and $\hat{\mathbf{n}}$ is the component normal to the motion, in the orbital plane), as a function of the time-to-impact $\Delta t$ expressed as a multiple of the NEO orbital period. The out of plane component $h$ of $\delta\hat{\mathbf{v}}_{opt}$ is not shown since it is always $<10^{-15}$.

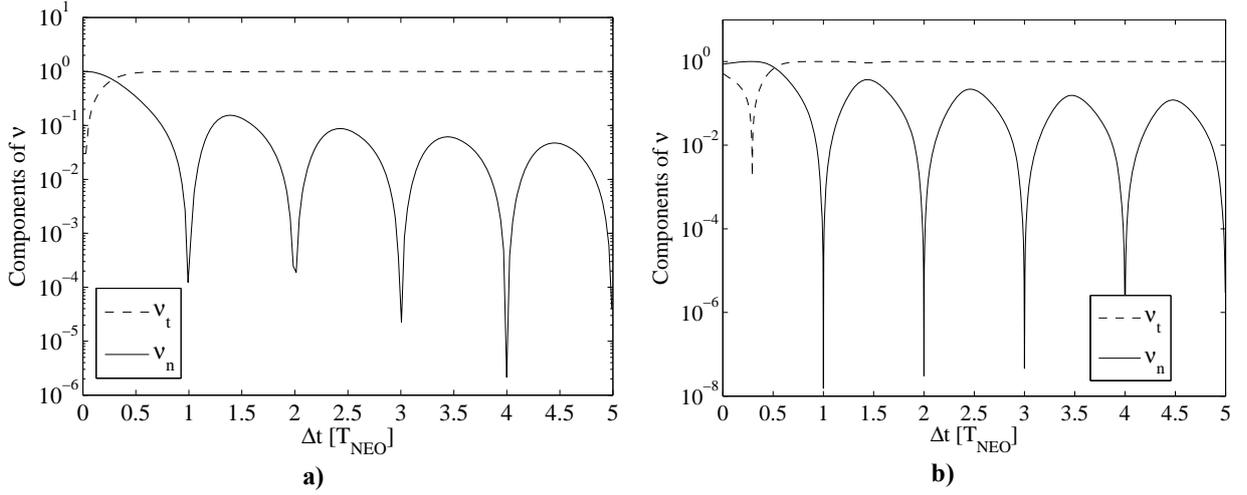

**Fig. 1 Components of the optimal δv direction for a) asteroid 2000SG344, b) asteroid 1979XB.**

As a result of this analysis we can derive that: for a $\Delta t$ smaller than a specific $\Delta t_{NEO} < 1 T_{NEO}$, different for every asteroid, the component perpendicular to the motion dominates the other two while for larger $\Delta t$, the tangential component becomes dominant. These conclusions are consistent with the results presented by Park and Ross[7] and with the ones of Conway[6]. It can be noted that the value of the normal component of the optimal deviation impulse goes periodically to zero with a period equal to the one of the asteroid. Therefore a deviation impulse given in the normal direction yields no deviation after an exact number of revolutions. Fig. 2a and b emphasize the optimality of the solution: the deviation obtained with $\|\delta\mathbf{v}\| = 0.07$ m/s was calculated, applying the maneuver along the optimal direction (solid line), and along the tangent (dot line), the normal (dash-dotted line) and the out-of-plane direction (dashed line). $\delta r$ associated to $\delta\hat{\mathbf{v}}_{opt}$ is the maximum one and overlaps the deviation achieved with a normal impulse, for low $\Delta t$, and the one with a tangent maneuver, for longer $\Delta t$. An impulsive action at the pericenter is found to be the most effective one, while a $\delta\mathbf{v}$ at the apocenter gives a deviation close to the minimum one. The choice of an optimum timing along the orbital period is more significant for highly eccentric orbits (see Fig. 2b).

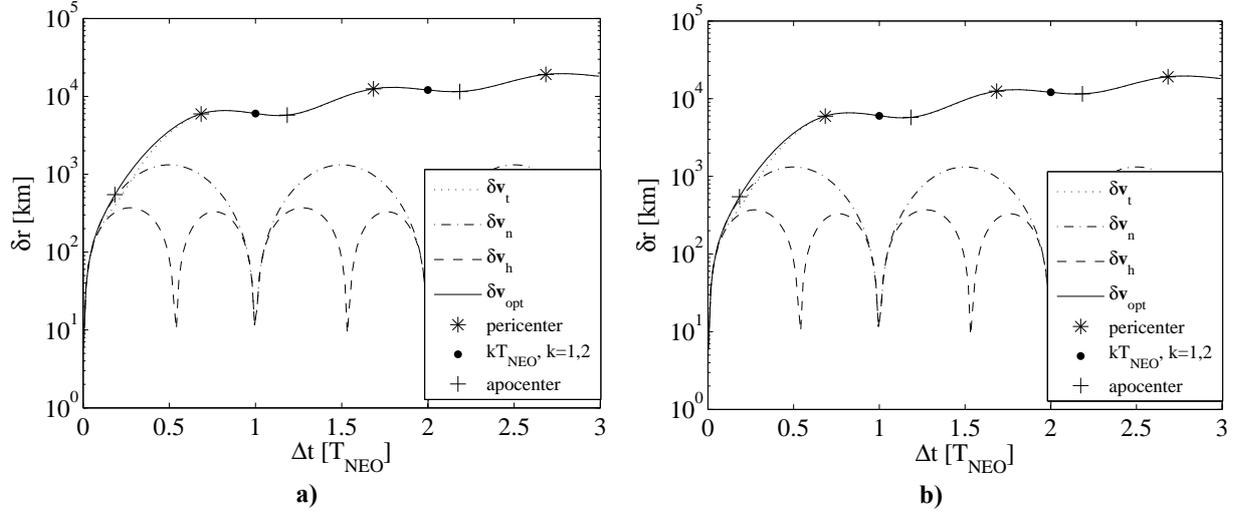

**Fig. 2  Deviation achieved with ||δv||=0.07 m/s for a) asteroid 2000SG344, b) asteroid 1979XB.**

While the direction of the optimal impulse is given by the maximization of the quadratic form associated with the transition matrix, the sign of $\delta\mathbf{v}_{opt}$, at this point of the analysis, is completely arbitrary. However, if we define the relative difference $e_{r1}$ between the deviation computed for $\delta\hat{\mathbf{v}}_{opt}$ and for $-\delta\hat{\mathbf{v}}_{opt}$, as,

$$e_{r1} = \frac{\left\|\delta\mathbf{r}_{+\hat{\delta}\mathbf{v}_{opt}} - \delta\mathbf{r}_{-\hat{\delta}\mathbf{v}_{opt}}\right\|}{\left\|\delta\mathbf{r}_{-\hat{\delta}\mathbf{v}_{opt}}\right\|}.$$

and we plot it as a function of the time-to-impact (see Fig. 3) we can conclude that the sign of $\delta\mathbf{v}_{opt}$ does not change the magnitude of the deviation. This can be alternatively demonstrated by changing the sign of $\delta\mathbf{v}$ in Eqs. (2). The variation of the orbital parameters is of opposite sign and consequently the displacement of the asteroid, described by Eqs. (1), is also in the opposite direction but with the same magnitude. This result confirms the results obtained by Conway[6].

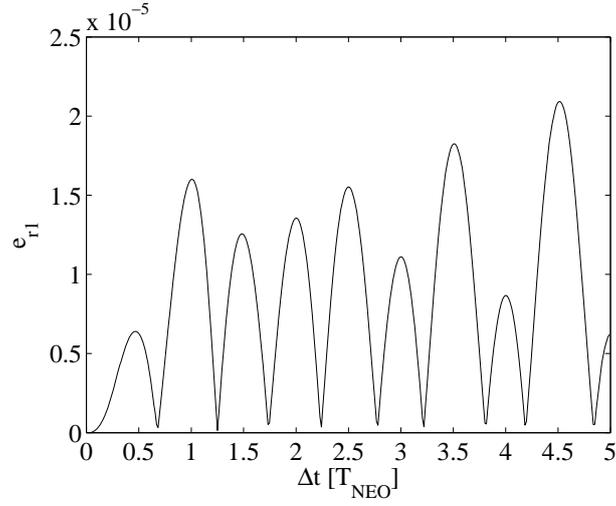

**Fig. 3 Relative error calculated by orbit propagation for 2000SG344.**

## B. Accuracy Analysis

The accuracy of Eqs. (1) was assessed by numerically propagating forward in time the deviated orbit of the asteroid and comparing the obtained variation in the position vector with the one predicted by Eqs. (1). The nominal trajectory was propagated from the deviation point up to the MOID, for a period up to 15 years, and the deviated one was integrated starting from the deviation point on the asteroid orbit with the perturbed velocity vector $\mathbf{v} + \delta\mathbf{v}$. As a measure of accuracy we computed the relative error between the variation in position after numerical propagation and the analytically estimated one:

$$e_r = \frac{\left\| \delta\mathbf{r}_{propagated} - \delta\mathbf{r}_{estimated} \right\|}{\left\| \delta\mathbf{r}_{propagated} \right\|} \tag{10}$$

Fig. 4 shows the relative errors, for the asteroid 2000SG344 and for the asteroid 1979XB, as a function of the time-to-impact $\Delta t$ and $\delta\mathbf{v}_{opt}$. These two asteroids, the former with $e < 0.1$ and $i < 10\,\text{deg}$, the latter with $e > 0.1$ and $i > 10\,\text{deg}$, were chosen in order to study the impact of the orbital parameters on the relative error.

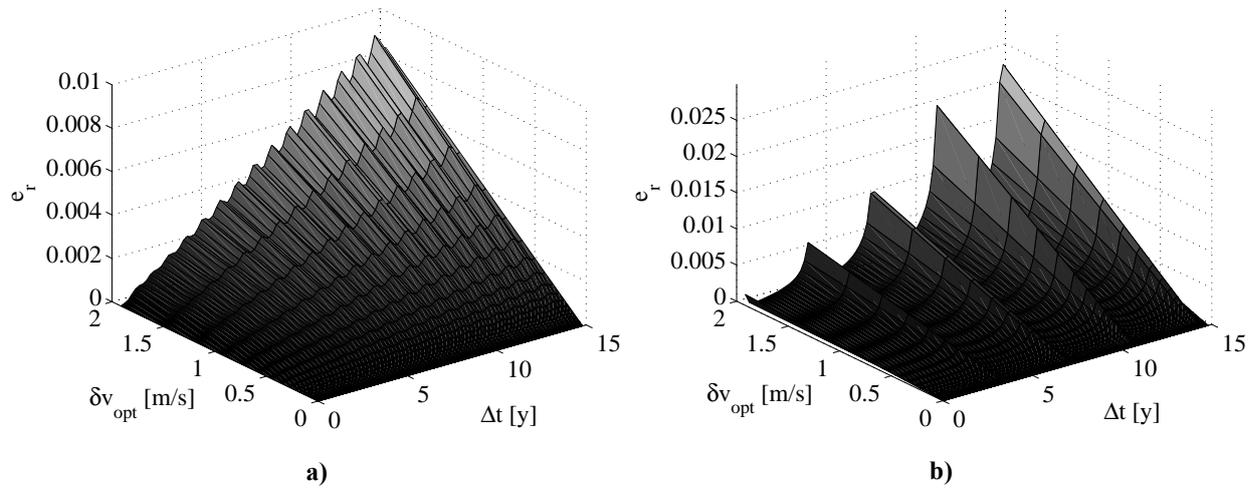

**Fig. 4 Relative error for the deviation of a) asteroid 2000SG344, b) asteroid 1979XB.**

For both asteroids the relative error grows with the time-to-impact and with $\|\delta \mathbf{v}_{opt}\|$ since the difference between the deviated orbit and the nominal one grows significantly and the proximal motion equations become inaccurate in describing the actual motion of the asteroid. In fact, as also stated by Schaub et al.[2], the hypotheses under which the equations were derived hold true until the relative orbit radius is small compared to the chief orbit radius. It is remarkable the difference in the maximum relative error between the two asteroids. If we compute the maximum relative error, i.e. for $\Delta t \leq 15$ y and $\delta v = 2$ m/s, for a large number of asteroids characterized by different sets of orbital parameters, we can see (Fig. 5) that its value grows with the eccentricity of the orbit of the asteroid.

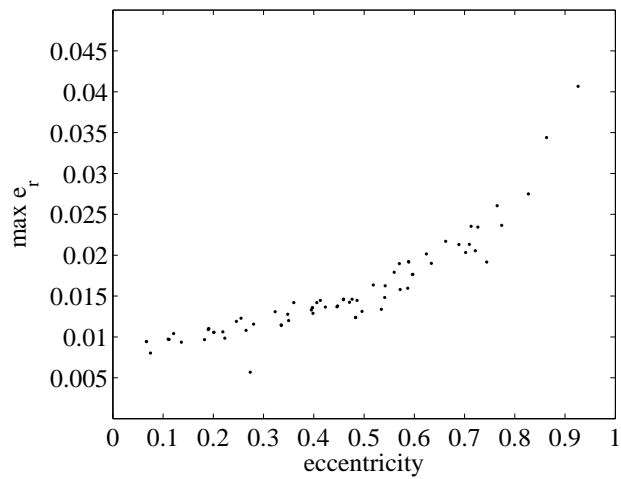

**Fig. 5 Maximum relative error for different asteroids.**

## C. Representation on the b-plane

The achieved deviation can be represented on the plane perpendicular to the incoming relative velocity of the small body at the planet arrival, i.e. the b-plane. We can define a local reference system centered in the Earth with the axis $\eta$ perpendicular to the b-plane aligned along the unperturbed velocity of the asteroid relative to the Earth, the axis $\zeta$ along the direction opposite to the projection of the heliocentric velocity of the planet onto the b-plane, and the axis $\xi$ that completes the reference system (see Fig. 6a). The general transformation from the Cartesian to the b-plane reference frame is:

$$\mathbf{x}_{b-plane} = \begin{bmatrix} \hat{\xi} & \hat{\eta} & \hat{\zeta} \end{bmatrix}^T \mathbf{x}_{cartesian}$$

where $\mathbf{x}$ is a generic vector, $\hat{\eta}$, $\hat{\xi}$ and $\hat{\zeta}$ are column vectors that can be computed as:

$$\hat{\eta} = \frac{\mathbf{U}_{NEO,nominal}}{\|\mathbf{U}_{NEO,nominal}\|}$$

$$\hat{\xi} = \frac{\mathbf{v}_E \times \hat{\eta}}{\|\mathbf{v}_E \times \hat{\eta}\|}$$

$$\hat{\zeta} = \hat{\xi} \times \hat{\eta}$$

where $\mathbf{U}_{NEO,nominal}$ is the unperturbed velocity of the asteroid relative to the Earth, expressed in a Cartesian reference frame and $\mathbf{v}_E$ is the velocity of the Earth.

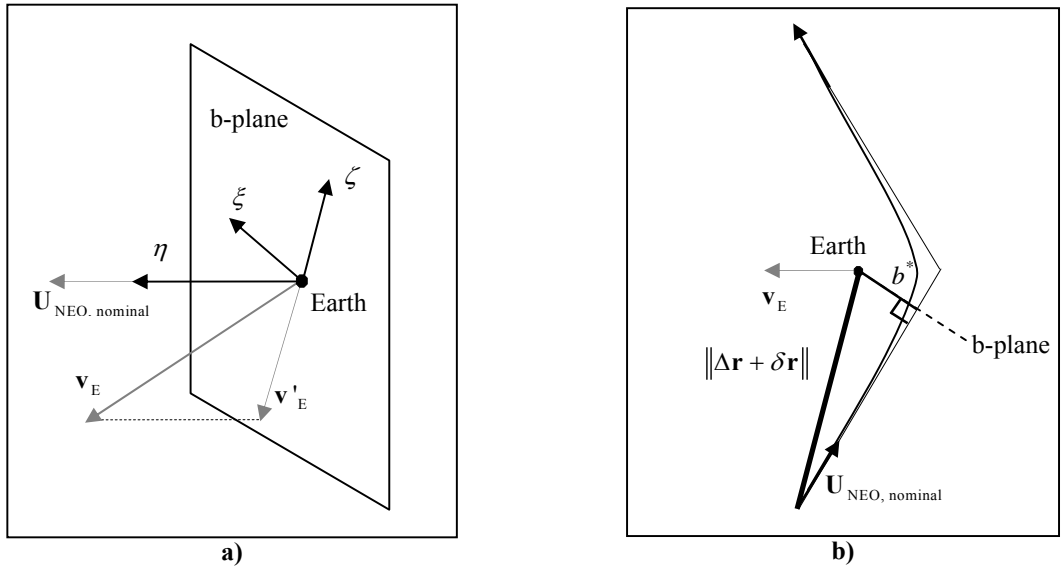

**Fig. 6 a) Representation on the b-plane. b) Geometry of hyperbolic passage.**

The proper representation would be on the *instantaneous* b-plane, perpendicular to the *deviated* relative velocity of the asteroid, however the maximum relative error between the plane perpendicular to the nominal relative velocity and the perturbed one is around 0.01, thus in the following we are using the former one which avoids the additional calculation of the velocity of the deflected asteroid. Moreover, for this analysis, we set to zero the distance at the MOID in order to have the Earth at the origin of the reference system on the b-plane. To this aim, the phase and periapsis anomaly of the asteroids were modified in order to have $\Delta r = 0$. This will not change the result of our analysis, because the other geometric properties of the orbit are unchanged.

On the b-plane we can represent the distance $b^*$ (called impact parameter) from the Earth to the intercept of the asymptote of the hyperbola of the deviated orbit of the asteroid:

$$b^* = \sqrt{\xi^2 + \zeta^2}$$

Since we expect the hyperbolic trajectory to be close to a straight line, $b^*$ is also a good estimation of the radius of the pericenter of the hyperbolic trajectory and hence the minimum distance from the Earth (see Fig. 6b). Fig. 7 shows the impact parameter (bold lines) for a highly elliptic asteroid (1979XB) together with the norm of the deviation δ**r** (thin lines). Although for a time-to-impact $\Delta t$ above a specific value, which is different for every asteroid (in the case of 1979XB, $\Delta t_{NEO} > 0.25 T_{NEO}$), the maximization of the $b^*$-parameter and the deviation lead to the same conclusion on the optimal deflection strategy. For smaller $\Delta t$ the b-plane suggests a different strategy. This can be appreciated in Fig. 8a that contains a close up of Fig. 7 for $\Delta t < 1 T_{NEO}$.

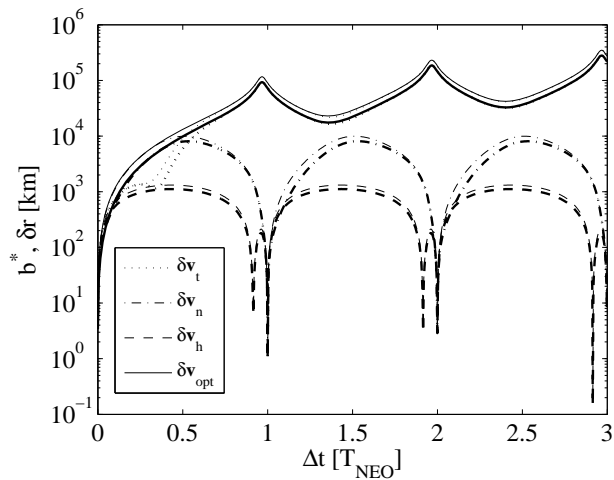

**Fig. 7 Impact parameter and magnitude of the deviation for 1979XB with δ$v$ = 0.07 m/s. The bold lines show the $b^*$-parameter, the thin lines the deviation.**

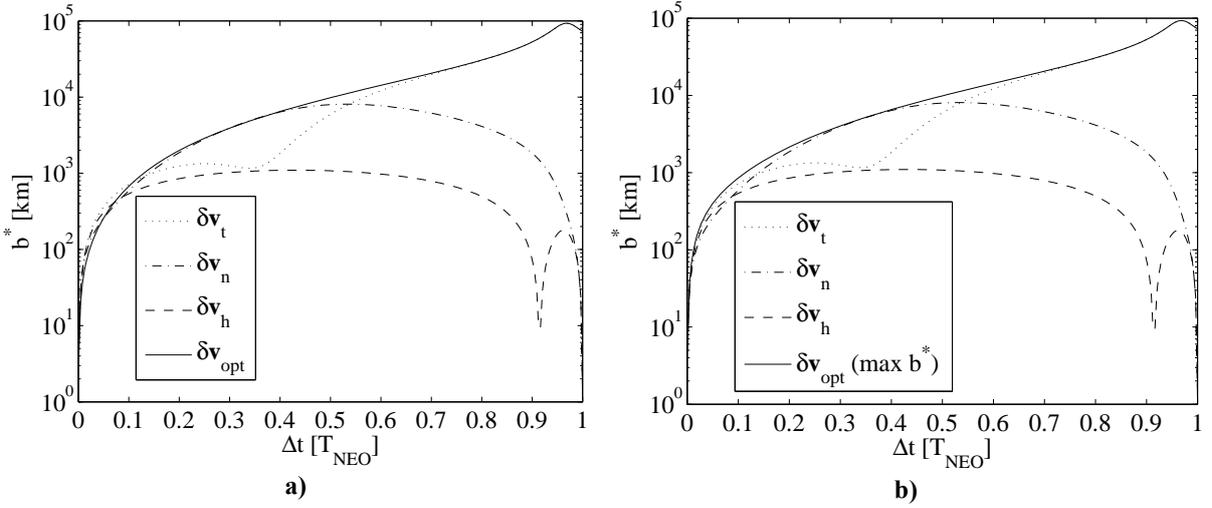

**Fig. 8 Impact parameter for 1979XB $\Delta t < 1 T_{NEO}$.
a) Strategy of maximum deviation. b) Strategy of maximum $b^*$-parameter.**

The difference between the two results is mainly due to the gravitational effect of the Earth when the asteroid approaches and enters the sphere of influence. The formulation of the maximization problem in Eq. (7) is modified in order to maximize the projection of the deviation in the b-plane; the result can be seen in Fig. 8b. For example for asteroid 1979XB, from the b-plane analysis we can conclude that the direction of the optimal impulse changes from the tangent direction, for a very small $\Delta t$, to the normal to the motion, for $\Delta t \cong 0.3 T_{NEO}$, while for $\Delta t > 0.7 T_{NEO}$ the tangential component becomes the dominant one.

In Fig. 9 the result for the $\delta r$-maximization strategy (thin line) is compared to the ones for the $b^*$-maximization strategy (bold line), for 2000SG344 in the range $\Delta t > 0.5 T_{NEO}$. The maximization of the impact parameter would lead to choose the *h*-direction strategy for very small $\Delta t$, the *n*-direction for a range of $0.15 T_{NEO} < \Delta t < 0.25 T_{NEO}$ and the tangential direction for higher $\Delta t$. Note that for small $\Delta t$ the angle of the optimal impulse changes depending on the orbital parameters of the asteroid, but for higher $\Delta t$ the optimal strategy is always the one along the direction of motion.

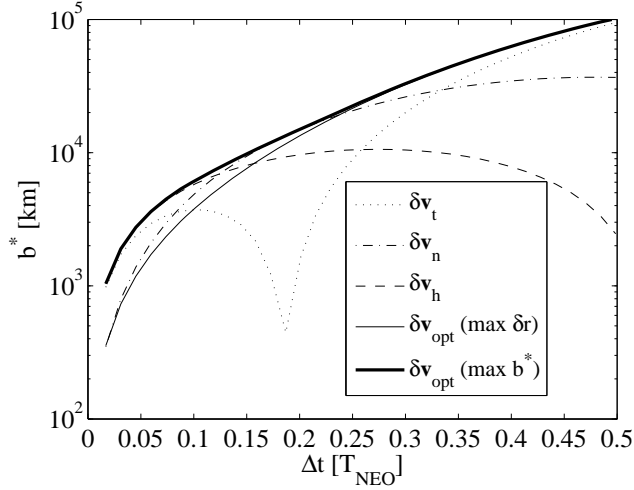

**Fig. 9 Impact parameter for 2000SG344 Δt<0.5T$_{NEO}$.
Strategy of maximum deviation (solid line) and maximum $b^*$-parameter (bold line).**

*1. Three-body analysis*

The results obtained with the b-plane formulation imply an increase of the $\delta\mathbf{v}$-requirement due to the gravitational effects of the Earth, which is consistent with the results found by Park and Ross[8]. Furthermore they suggest a different optimal strategy for short times-to-impact. We can verify the reliability of these results by propagating the motion of the asteroid after the deflection maneuver for two different cases: the optimal deflection maneuver is computed as the result of the maximization of the deviation, the optimal deflection maneuver is computed as the result of the maximization of the $b^*$-parameter. A full three-body dynamic model was used (i.e. considering the Sun and the Earth as gravitational bodies) and the trajectory was propagated, after a 2 m/s impulse, over the interval $I_p = \begin{bmatrix} t_d & t_{MOID} + 0.1 T_{NEO} \end{bmatrix}$. Then, the closest point to the Earth was computed as follows:

$$\delta r_{\min,3b} = \min_{t \in I_p} \left\| r_E(t) - r_{NEO}(t) \right\| \tag{11}$$

As can be seen in Fig. 10a, a deflection maneuver computed maximizing the deviation is not an optimal strategy for short times-to-impact while the one computed maximizing the $b^*$-parameter (see Fig. 15b) leads to better results. Note that this is true for short times-to-impact while for longer ones the two strategies are equivalent. With both strategies the projection of the deviation on the b-plane (bold line) is a reliable estimation of the actual deviation computed with the three-body model (thin line). The deviation at the MOID, considering the two-body dynamics (dashed line), instead, does not predict accurately the actual minimum distance from the asteroid in proximity of the Earth.

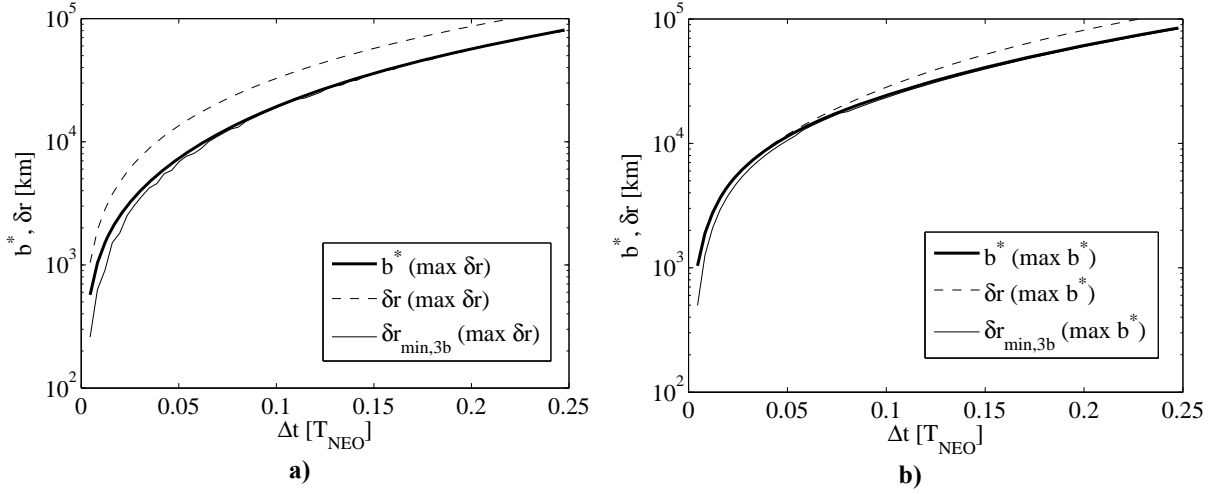

**Fig. 10** a) Maximum deviation strategy for 1979XB, b) Maximum $b^*$-parameter strategy for 1979XB. The dashed line and the bold line represent respectively the deviation and its projection on the b-plane, calculated through the two-body problem. The continuous thin line represents the minimum deviation computed through the three-body problem.

*2. Analysis of the Deviation Components in the b-plane*

Fig. 11 shows the components of the deviation in the b-plane as a function of the time-to-impact, when the optimal strategy is computed maximizing $b^*$. In the same figure $t_{\delta r_{min,3b}}$ is the time corresponding to $\delta r_{min,3b}$, defined in Eq. (11); $t_{\delta r_{min,3b}} - t_{MOID}$ represents the difference between the instant when the actual minimum distance from the Earth is reached and the expected time at the MOID. This quantity is expressed in days, multiplied by $10^6$ to make it comparable in scale with the components of the deviation.

The components of the deviation projected onto the b-plane have a discontinuity corresponding to the discontinuity in $t_{\delta r_{min,3b}} - t_{MOID}$. In particular when $t_{\delta r_{min,3b}} - t_{MOID} > 0$, we have $\eta < 0$; this means that the asteroid at $t_{MOID}$ has not intersected the b-plane yet (the component normal to it is negative). This situation is depicted in Fig. 12a, where the point A represents the asteroid approaching a fly-by of the Earth. When $t_{\delta r_{min,3b}} - t_{MOID} < 0$, instead, then $\eta > 0$; this means that the asteroid at $t_{MOID}$ has already intersected the b-plane (the component normal to it is positive). This situation is depicted in Fig. 12b, where the point B represents the asteroid after the fly-by.

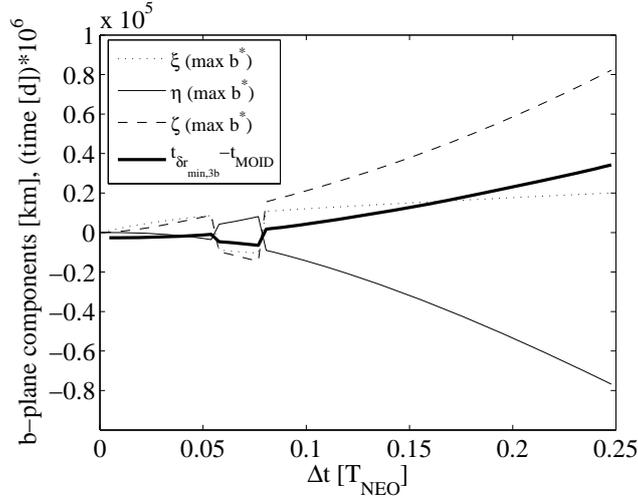

**Fig. 11    Components of the deviation in the b-plane for asteroid 1979XB.**

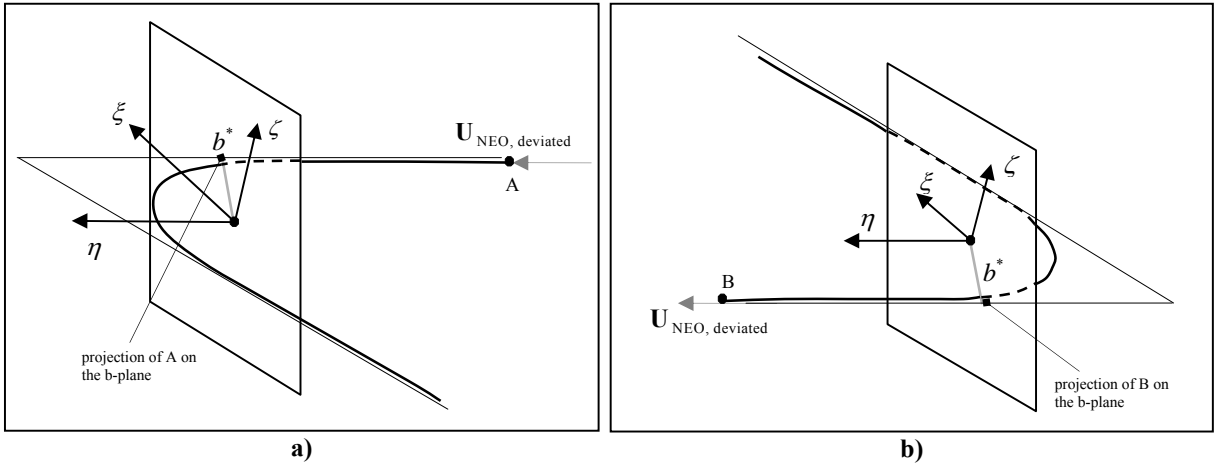

**Fig. 12    Fly-by representation in the b-plane reference system. a) Case A: the asteroid is approaching a fly-by of Earth. b) Case B: the asteroid at the end of the fly-by.**

The two situations described in Fig. 12 are a consequence of the sign of the impulsive maneuver $\delta \mathbf{v}_{opt}$ obtained from problem (7). In fact, for a $+\delta \mathbf{v}$ along the motion, the period of the asteroid is increased, hence at $t_{MOID}$ the asteroid is at point A in Fig. 12a. On the other hand for a $-\delta \mathbf{v}$ along the motion, the period of the asteroid is decreased, hence at $t_{MOID}$ the asteroid is at point B in Fig. 12b.

Fig. 13 represents the projection on the b-plane of the deviated points for different values of $\Delta t$. The deviation was calculated by applying the impulsive maneuver along the tangent (Fig. 13a), the normal (Fig. 13b) and the out–of-plane direction (Fig. 13c) respectively. It can be noted that an impulse along the tangent direction produces a substantial variation of the $\zeta$-component, with a secular and a periodic term, and a small, periodic variation of the

$\xi$-component. An impulse in the normal direction, instead, produces a purely periodic variation of both components. The difference in the three strategies can be appreciated in Fig. 14a (asteroid 1979XB) and Fig. 14b (asteroid 20000SG344), that represent the development of the three components along the time axis.

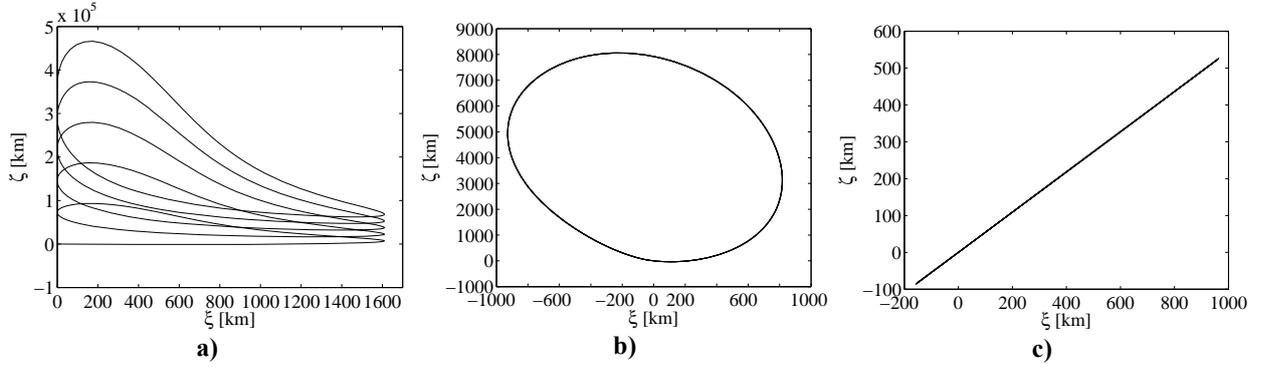

**Fig. 13   Projection on the b-plane of the deviation for 1979XB.**
**δv = 0.07 m/ s applied a) along the tangent to the motion, b) along the normal to the motion, c) along the *h*-direction.**

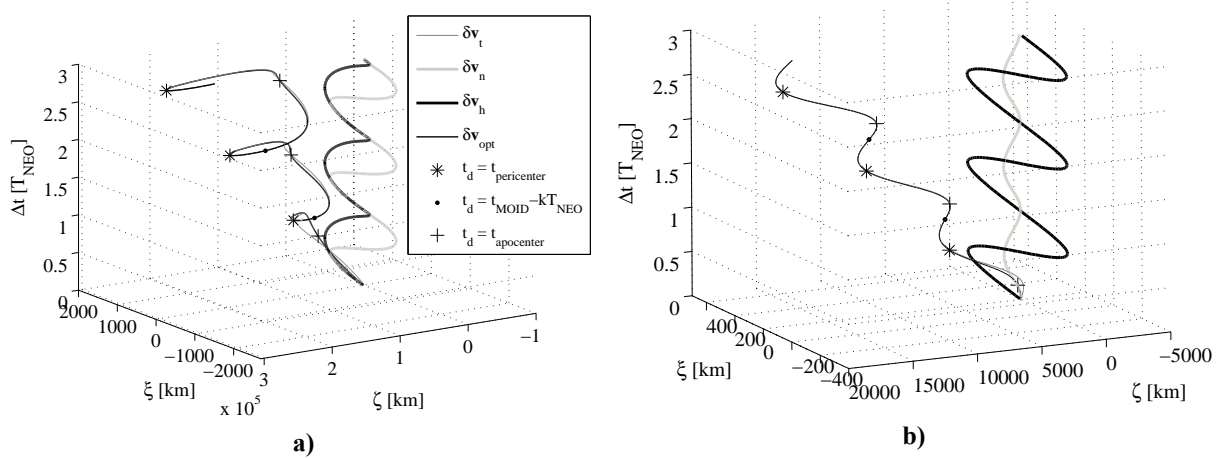

**Fig. 14   Projection on the b-plane, function of Δt for a) asteroid 1979XB, b) asteroid 2000SG344.**

Fig. 15 shows on the same plot the result of a deflection action given in the three directions for a highly elliptic orbit (asteroid 1979XB) and low-eccentric orbit (asteroid 20000SG344). It can be noted that, for high elliptical orbit, such as 1979XB one (see Fig. 15a), the best results are achieved if the impulse is given at the pericenter of the orbit. A $\delta \mathbf{v}$ at the apocenter, instead, is almost the less efficient action, because it changes $\xi$ (related to the MOID) but not $\zeta$. Note that by acting *k* orbital periods before the impact (where *k* is an integer), the deviation component along $\xi$ is zero. For a orbit with a low eccentricity (see Fig. 15b) a deviation maneuver at the pericenter is still the most efficient, though it does not maximize $\zeta$.

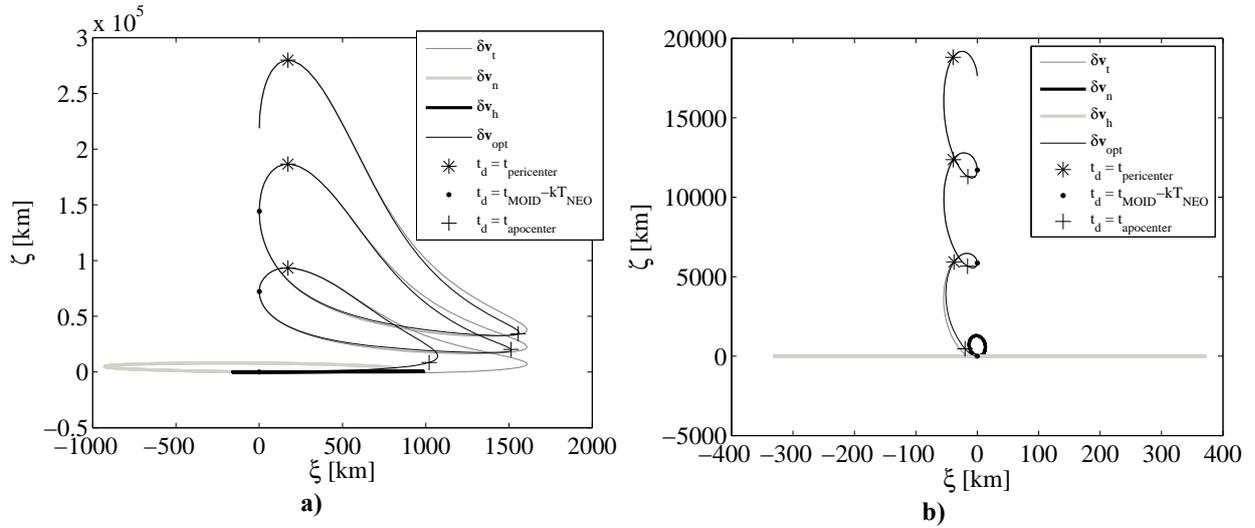

**Fig. 15** Projection on the b-plane of the deviation. δv = 0.07 m/s applied along the optimal direction (normal line) for a) asteroid 1979XB, b) asteroid 2000SG344.

## III. Mission Options

### A. Target Selection

Potentially Hazardous Asteroids (PHAs) are a subclass of NEOs which are defined based on parameters that measure the asteroids' potential to make dangerously close approaches to Earth. Specifically all asteroids having a MOID smaller than 0.05 AU and an absolute magnitude of 22 or less (i.e. bigger than about 150 m in diameters) define this group. More than 770 PHAs are currently known, and it is estimated that at least a similar quantity of objects with these characteristics are still to be discovered. Different research groups in the world, such as the JPL Sentry system[‡] or the NeoDys group in Pisa Italy[§], keep updated databases continuously assessing the risk posed by these objects. As more ground based observations become available a more accurate determination of the PHAs orbits will be available and, as a consequence, some of the asteroids might be removed from the possible impact risks list. However looking at the current estimation for potential impacts through the next century, there are still a few objects whose impact probability is not negligible in statistical terms.

Table 1 shows an extract of 30 objects taken from the JPL catalogue of asteroids, with the exception of 2004VD17 that, although foreseen to have its closest approach with the Earth in the years 2102-2104, is currently considered as the most dangerous objects in terms of impact probability due to its large size relative to the other PHAs listed.

---

[‡] http://neo.jpl.nasa.gov/risk/
[§] http://newton.dm.unipi.it/cgi-bin/neodys/neoibo

Table 1 lists the asteroids, and their properties, used in the analysis. For simplicity, each of the asteroids is given a local reference number (instead of using the formal names or international ID numbers). The semi-major axis $a$ is given in astronomical units and the inclination $i$, anomaly of the pericenter $\omega$, argument of the ascending node $\Omega$ and mean motion $M$ are in degrees, the estimated mass in kg and the epoch is in Modified Julian Days (MJD). The list considered contains some bodies that have recently become objects of interest for the scientific community: Apophis, 2004VD17 and 2005WY5 are in fact reported in the JPL catalogue of the most recently observed objects as the currently most risky, having a Palermo scale ranging between –0.57 and –2.61. Some of the objects in the list are instead among those not recently observed or even lost. This is in fact a major issue in the current assessment of their risk; due to the limited capabilities of ground based observation and to limited available resources, most of the hazardous objects can be lost for even many years, this resulting in the possibility that when new observations of the objects are available again, these could definitely rule out the possibility of an impact or actually turn out to have an increased impact probability, with the additional drawback of a reduced warning time.

The rationale behind the selection presented in Table 1 is twofold. We are interested in providing some general considerations on optimal impact trajectories and consequent deviations strategies, and for this reason have surveyed a set of potential dangerous objects presenting a large variety of orbital and physical characteristics. As can be noticed from Table 1, our selection collects objects having semi-major axes ranging between 0.85 and 3 AU, an eccentricity as large as 0.92 and an orbital inclination up to 28 deg, while having estimated masses in a range between $10^7$ and $10^{12}$ kg. Such relevant differences in both orbital elements and the mass will eventually affect the optimal impact and deflection strategy. At the same time we want to look at some actual sample cases, considering real objects currently having quite a high impact probability (usually indicated in terms of Palermo scale), in order to provide a worst case assessment of the current and short term future capabilities of deflecting hazardous objects like Apophis or 2004VD17 if such an unlikely, but highly disastrous event, should ever be faced by our society. The MOID $\Delta r$ was calculated using the Earth's ephemerides on the 1st of January 2000, at 12:00 (0 MJD since 2000). As a consequence of this approximation, the MOID of asteroid 1997XR2 results to be less than the Earth radius. The actual MOID varies with time, due to the actual orbit of both the Earth and the asteroid, furthermore a MOID smaller than the radius of the Earth does not imply an imminent impact since the Earth and the asteroid could not be at the MOID at the same time, Note that, the aim of this work is not to reproduce a realistic impact scenario but rather to assess the actual achievable deviation, as opposed to the theoretical one derived in the previous chapters, depending on the mass and orbital characteristics of the asteroid. In this respect the modulus and direction of the MOID vector play an important role, as it will be demonstrated in the following. A more accurate calculation of the MOID would produce a more precise estimation of the actual achievable deviation but would not invalidate the results in this paper.

Table 1    Physical parameters for considered NEOs.

| Id | Name | $a$ | $e$ | $i$ | $\omega$ | $\Omega$ | $M$ | epoch [MJD] | mass [kg] | $\Delta r$ [km] |
|---|---|---|---|---|---|---|---|---|---|---|
| 1 | 2004VD17 | 1.50 | 0.58 | 4.22 | 90.7 | 224.2 | 286.9 | 53800.5 | 2.7e11 | 229479.20 |
| 2 | Apophis | 0.92 | 0.19 | 3.33 | 126.3 | 204.4 | 222.2 | 53800.5 | 4.6e10 | 36651.75 |
| 3 | 2005WY55 | 2.47 | 0.72 | 7.26 | 285.9 | 248.4 | 3.30 | 53800.5 | 1.9e10 | 696520.60 |
| 4 | 1997XR2 | 1.07 | 0.20 | 7.17 | 84.6 | 250.8 | 211.8 | 53800.5 | 1.7e10 | 3277.43 |
| 5 | 1994WR12 | 0.75 | 0.39 | 6.81 | 205.8 | 62.8 | 27.3 | 53700 | 2.0e9 | 283313.30 |
| 6 | 1979XB | 2.35 | 0.72 | 25.1 | 75.7 | 85.5 | 62.0 | 53700 | 4.4e11 | 3720840.42 |
| 7 | 2000SG344 | 0.97 | 0.06 | 0.11 | 274.9 | 192.3 | 132.3 | 53800.5 | 7.1e7 | 124351.73 |
| 8 | 2000QS7 | 2.68 | 0.66 | 3.19 | 218.7 | 153.5 | 84.8 | 53800.5 | 9.9e10 | 542496.18 |
| 9 | 1998HJ3 | 1.98 | 0.74 | 6.54 | 92.7 | 224.9 | 333.6 | 50926.5 | 4.5e11 | 1907030.74 |
| 10 | 2005TU45 | 1.97 | 0.49 | 28.5 | 76.8 | 120.2 | 34.1 | 53651.5 | 3.3e12 | 38152163.70 |
| 11 | 2004XK3 | 1.21 | 0.25 | 1.43 | 302.2 | 58.1 | 22.0 | 53800.5 | 1.1e8 | 168758.33 |
| 12 | 1994GK | 1.92 | 0.59 | 5.60 | 111.4 | 15.4 | 17.3 | 49450.5 | 1.5e8 | 445443.47 |
| 13 | 2000SB45 | 1.55 | 0.39 | 3.67 | 216.3 | 195.5 | 214.4 | 53700 | 1.3e8 | 199226.54 |
| 14 | 2001CA21 | 1.66 | 0.77 | 4.93 | 218.8 | 46.4 | 65.5 | 53700 | 4.3e11 | 5574409.52 |
| 15 | 2005QK76 | 1.40 | 0.51 | 22.9 | 266.1 | 337.6 | 36.1 | 53613.5 | 4.1e7 | 122907.28 |
| 16 | 2002TX55 | 2.23 | 0.57 | 4.37 | 148.8 | 190.2 | 16.8 | 53800.5 | 3.4e8 | 534543.89 |
| 17 | 2005EL70 | 2.27 | 0.92 | 16.18 | 167.5 | 167.5 | 12.0 | 53438.5 | 1.9e8 | 21308100.31 |
| 18 | 2001BB16 | 0.85 | 0.17 | 2.02 | 195.5 | 122.5 | 327.4 | 53800.5 | 1.5e9 | 704667.59 |
| 19 | 2002VU17 | 2.47 | 0.61 | 1.49 | 308.75 | 55.67 | 11.37 | 52599.5 | 7.3e7 | 1500966.15 |
| 20 | 2000TU28 | 1.07 | 0.18 | 15.64 | 280.6 | 203.1 | 227.0 | 53800.5 | 3.0e10 | 166332.26 |
| 21 | 2001AV43 | 1.27 | 0.23 | 0.27 | 43.0 | 30.7 | 226.9 | 53800.5 | 1.2e8 | 632550.85 |
| 22 | 2002RB182 | 2.54 | 0.65 | 0.22 | 254.3 | 165.5 | 347.4 | 52532.5 | 1.1e9 | 302338.44 |
| 23 | 2002GJ8 | 2.97 | 0.82 | 5.30 | 180.3 | 144.2 | 261.3 | 53800.5 | 1.3e11 | 13925769.75 |
| 24 | 2001FB90 | 2.48 | 0.78 | 1.92 | 14.5 | 266.3 | 343.3 | 51993.5 | 5.7e10 | 4781828.50 |
| 25 | 2005NX55 | 1.52 | 0.58 | 26.16 | 277.2 | 106.4 | 327.2 | 53563.5 | 3.8e9 | 5098118.30 |
| 26 | 1996TC1 | 1.86 | 0.72 | 14.53 | 258.8 | 5.01 | 22.8 | 50363.5 | 2.3e8 | 11305879.51 |
| 27 | 6344P-L | 2.64 | 0.64 | 4.66 | 232.6 | 184.9 | 349.8 | 37203.5 | 1.2e10 | 4183900.25 |
| 28 | 2004ME6 | 2.36 | 0.57 | 9.44 | 210.3 | 112.2 | 346.1 | 53182.5 | 1.5e9 | 4343813.94 |
| 29 | 2001QJ96 | 1.59 | 0.79 | 5.87 | 121.3 | 339.1 | 333.9 | 52147.5 | 3.3e9 | 292749.39 |
| 30 | 2004GE2 | 2.04 | 0.70 | 2.16 | 259.9 | 45.1 | 341.6 | 53112.5 | 8.0e9 | 856426.32 |

**B. Impact Model and Optimization Problem Definition**

The impact between the spacecraft and the asteroid is considered to be perfectly inelastic; we do not take into consideration other effects due to the ejection of mass or gasses. The variation of velocity imparted by the spacecraft to the asteroid is therefore given by the equation:

$$\delta \mathbf{v}_A = \gamma \frac{m_s}{m_A} \Delta \mathbf{v}_s \qquad (12)$$

where the relative velocity $\Delta \mathbf{v}_s$ of the spacecraft with respect to the asteroid at the impact point is computed from the ephemerides of the asteroid and from the solution of a Lambert's problem for the spacecraft, and the parameter $\gamma$ has a value of 1 in this implementation.

The mass of the asteroid $m_A$ was estimated from its measured magnitude, while the mass of the spacecraft at the impact point was computed through the rocket equation as follows:

$$m_s(t_d) = m_s(t_0) e^{-\frac{\Delta v_{tot}}{g_0 I_{sp}}} \qquad (13)$$

where the $I_{sp}$ was taken equal to 315 s and the total $\Delta v_{tot}$ is the sum of all the required maneuvers that the spacecraft has to perform after launch. All the celestial bodies are considered to be point masses with no gravity, the ephemerides of the asteroids were computed using the mean orbital elements in Table 1 while analytical ephemerides considering the long term variation of the orbital elements were used for the Earth and for Venus. In the case of direct Earth-asteroid transfers, $\Delta v_{tot}$ is the required velocity change at the Earth to reach the asteroid, in the case of Earth-Venus-asteroid transfers, instead, $\Delta v_{tot}$ accounts for the required velocity change at the Earth to reach Venus plus the deep space correction required after the Venus swing-by to reach the asteroid (further details on the trajectory model can be found in Vasile et al.[9]). The initial mass of the spacecraft is $m_s(t_0)$ = 1000 kg and the launcher is assumed to provide an escape velocity of 2.5 km/s. If the required $\Delta v_{tot}$ for the transfer is less than the escape velocity provided by the launcher, a higher effective mass at launch is considered, in order to fully exploit the launcher capabilities. In this case the initial mass is:

$$m_s(t_0) = 1000 e^{\frac{\Delta v_{exc}}{g_0 I_{sp}}}$$

where $\Delta v_{exc} = 2.5 - \Delta v_{tot}$. We consider that a minimum of 20% of the mass of the spacecraft at launch is allocated to structure and subsystems while a minimum of 10% of the propellant mass is allocated to tanks and propulsion system, therefore the quantity $1.1e^{-\frac{\Delta v_s}{g_0 I_{sp}}} - 0.3$ at impact must be positive. Hence we define a constraint $C_m$ on the residual mass computed at the impact:

$$C_m = \max\left(\left[sign\left(1.1e^{-\frac{\Delta v_s}{g_0 I_{sp}}} - 0.3\right), 0\right]\right) \qquad (14)$$

The deviation δ**r** is therefore a function of the mass of the spacecraft at impact and can be written in compact form as follows:

$$\delta\mathbf{r} = \gamma \frac{C_m m_s(t_d)}{m_A} \mathbf{T}\Delta\mathbf{v}_s \qquad (15)$$

The square of the modulus of the MOID after deviation is therefore:

$$J = \left(\Delta\mathbf{r} + \gamma \frac{C_m m_s(t_d)}{m_A} \mathbf{T}\Delta\mathbf{v}_s\right)^T \left(\Delta\mathbf{r} + \gamma \frac{C_m m_s(t_d)}{m_A} \mathbf{T}\Delta\mathbf{v}_s\right) \qquad (16)$$

which has to be maximized with respect to the launch date $t_0$ and the deviation time $t_d$. Note that, though from the analysis presented in the previous chapter the strategies that aims at maximizing $b^*$ are more accurate than the ones aiming at the maximization of δ*r*, in the following we use the latter since it provides good and reliable results for medium to long times-to-impact and requires a lower computational cost.

In order to better examine the full range of launch opportunities we ran three different optimizations fixing three different upper limits for the maximum warning time, which is the time from launch to the time the asteroid reaches the MOID: respectively up to 5, 10 and 15 years. This was obtained by fixing the upper limit for a possible launch date to respectively [3650 5475] MJD2000, [3650 7300] MJD2000 and [3650 9125] MJD2000. Since the warning time can be up to 15 years we computed all the times the asteroid is crossing the MOID for up to 15 years after the upper limit for the launch date and we took the first date the asteroid is reaching the MOID. Since some asteroids in Table 1 can not be reached with a low cost direct transfer, we also analyze the benefits of a single swing-by maneuver with Venus. More complex sequences and multi-burn maneuvers can further improve the deviation in the desired time frame; this will be the subject of future work.

Since we are interested in a large number of local minima for the objective function given in Eq. (16) rather than only the global one, we used a particular global optimization method which blends a stochastic search with an automatic solution space decomposition technique. This method has proven to be particularly effective compared to common optimization methods, especially when applied to space trajectory optimization problems[10,12].

Furthermore, it is expected that the largest deviations can be obtained with the longest warning time, thus we performed an additional analysis optimizing the warning time $t_w$ along with objective function (16):

$$t_w = t_{MOID} - t_0 \qquad (17)$$

For this second analysis the aim was to find the set of Pareto optimal solutions, i.e. all those solutions for which there is no other solution that has a better value for both $J$ and $t_w$. We used the same optimization method but in its multi-objective version[12,13].

### C. Results

*1. Single-objective optimization*

The results of the single-objective optimization consist of a number of families of mission opportunities for each upper boundary on the maximum warning time. When the asteroid has high inclination, the optimal interception points are concentrated at the ascending and descending node of the orbit. Two examples are shown in Fig. 16 (asteroid 1979XB and asteroid 1996TC1) with the value of the argument of latitude at interception shown in Fig. 17.

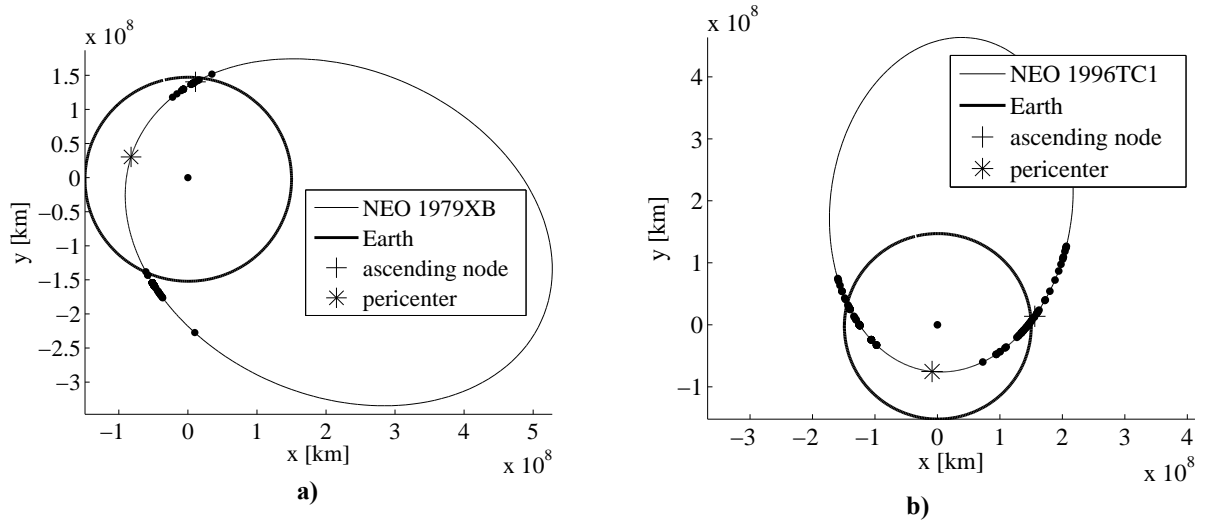

**Fig. 16   Optimal interception of a) asteroid 1979XB, b) asteroid 1996TC1.**

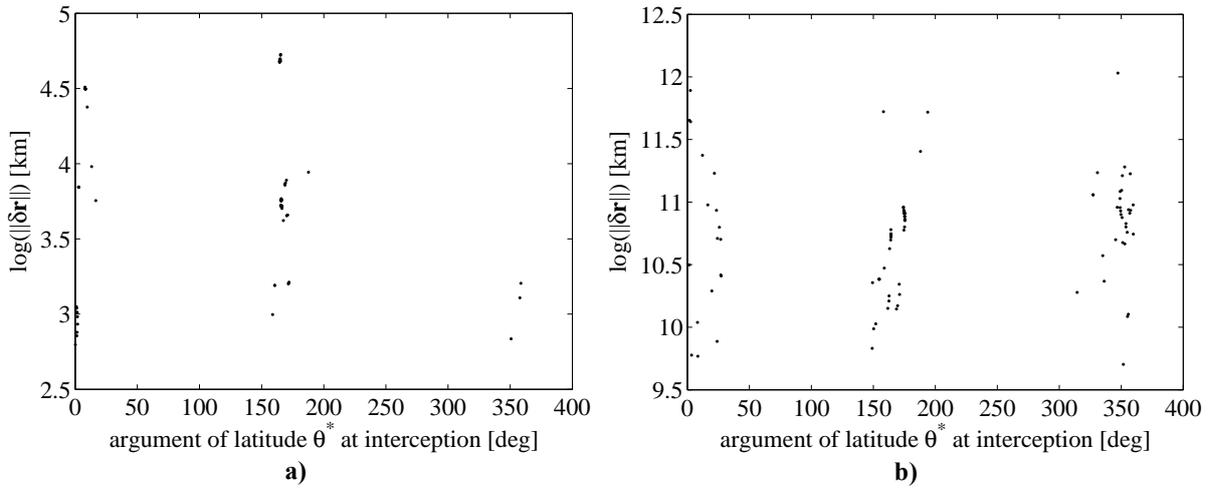

**Fig. 17** Optimal interception arguments of latitude for a) asteroid 1979XB, b) asteroid 1996TC1.

As can be seen in Fig. 16, the interception points, marked with a dot, are situated straddling the pericenter. This is the best compromise between an impact at the pericenter, which is the point that ensures the maximum change in the orbital period, and the transfer trajectory to reach the asteroid from the Earth. On the other hand, when the pericenter of the asteroid orbit is close to the Earth orbit, like in the case of asteroid 2000SB45 and 2002TX55 (see Fig. 18), many optimal solutions are grouped around the pericenter.

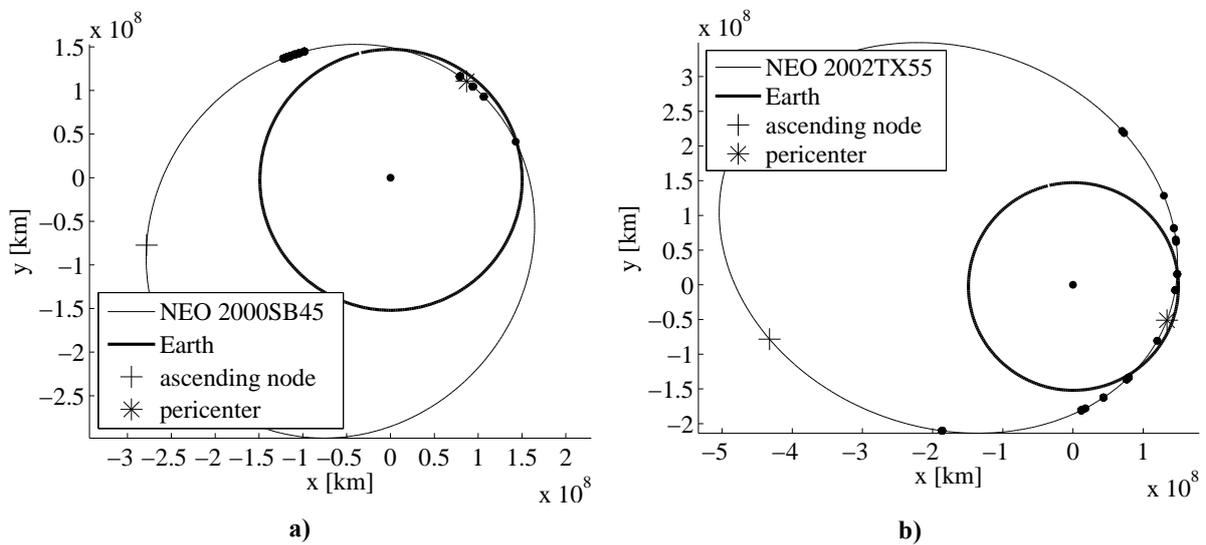

**Fig. 18** Optimal interception of a) asteroid 2000SB45, b) asteroid 2002TX55.

The value of the impact velocity is almost a linear function of the eccentricity of the orbit as shown in Fig. 19a and its out-of-plane component increases with the inclination of the orbit (see Fig. 19b).

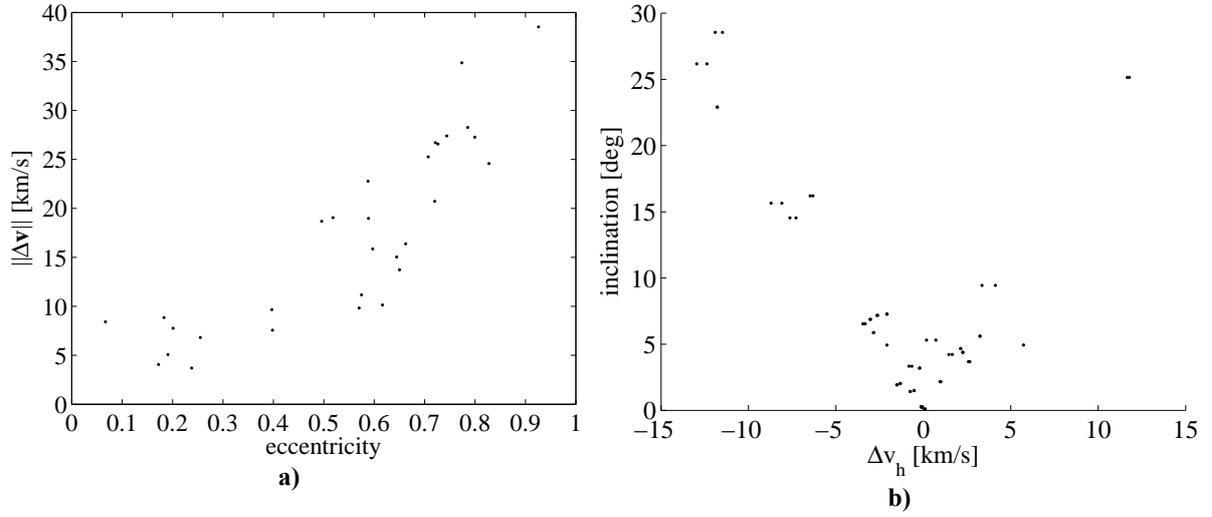

**Fig. 19    a) Impact velocity function of the eccentricity. b) *h*-component of the impact velocity function of the inclination.**

*2. Multi-objective optimization*

For each mitigation scenario, i.e. different upper boundary on the maximum warning time, a number of solutions were found that are Pareto optimal with respect to the total deviation $\|\Delta \mathbf{r} + \delta \mathbf{r}\|$ and warning time $t_w$. Two sets of Pareto optimal solutions are given as an example in Fig. 20. The former asteroid (2000SG344) has a low eccentric orbit, while the latter one (2002GJ8) has an eccentricity $e = 0.82$. Among all the solutions of each Pareto optimal set for each scenario, we selected, and listed in Table 2, the ones that maximize the total deflection. For each mission, details of the optimal trajectory are given, namely the launch date $t_0$, the time of flight $ToF$ and the mass of the spacecraft at the interception point with the asteroid $m_s(t_d)$. The close approach selected for each case is identified by $t_{MOID}$, that determines also the warning time for that mission $t_w$. $\Delta v_t$, $\Delta v_n$ and $\Delta v_h$ are the components of the relative velocity of the spacecraft with respect to the NEO at the interception point, while $\delta r$ is the achieved deviation.

**Fig. 20** Pareto front for a) asteroid 2000SG344, b) asteroid 2002GJ8.

**Table 2** Optimal launch opportunities for a direct transfer to selected asteroids. Result of the multi-objective optimization. The designation numbers correspond to the asteroids listed in Table 1.

| Id | $t_0$ [d] | ToF [d] | $t_{MOID}$ [MJD since 2000] | $t_w$ [d] | $m_s(t_d)$ [d] | $\Delta v_t$ [km/s] | $\Delta v_n$ [km/s] | $\Delta v_h$ [km/s] | $\delta r$ [km] | $\|\Delta\mathbf{r}+\delta\mathbf{r}\| - \|\Delta\mathbf{r}\|$ [km] |
|---|---|---|---|---|---|---|---|---|---|---|
| 1 | 4100.06 | 408.52 | 5826.16 | 1726.10 | 743.48 | -10.48 | -19.00 | 1.27 | 16.38 | 1.00 |
|   | 4836.37 | 341.55 | 7855.58 | 3019.20 | 962.04 | -10.05 | -17.83 | 1.25 | 42.24 | 2.20 |
|   | 4105.08 | 403.60 | 9208.52 | 5103.44 | 743.74 | -10.47 | -18.97 | 1.18 | 56.84 | 3.50 |
| 2 | 4165.11 | 310.87 | 5842.53 | 1677.42 | 1176.78 | 2.65 | -2.75 | -1.15 | 16.67 | 4.36 |
|   | 4697.60 | 62.89 | 7460.40 | 2762.80 | 940.69 | 3.24 | -0.78 | -0.82 | 35.82 | 10.73 |
|   | 4697.91 | 65.67 | 9401.85 | 4703.93 | 965.74 | 3.17 | -0.53 | -0.81 | 62.29 | 18.58 |
| 3 | 4671.63 | 507.57 | 9418.49 | 4746.86 | 481.39 | -16.60 | -20.91 | -2.15 | 1430.33 | 169.12 |
| 4 | 5448.36 | 146.33 | 5618.12 | 169.75 | 374.42 | -1.64 | 10.09 | 0.14 | 0.48 | 0.17 |
|   | 4946.31 | 245.89 | 7658.52 | 2712.20 | 527.53 | -3.63 | 6.75 | -2.63 | 75.76 | 12.21 |
|   | 5369.94 | 229.44 | 9290.84 | 3920.89 | 968.31 | -2.25 | 7.03 | -3.83 | 128.66 | 22.55 |
| 5 | 4103.70 | 334.67 | 5585.66 | 1481.96 | 943.91 | 6.66 | 1.81 | -2.94 | 489.45 | 122.24 |
|   | 4100.17 | 338.79 | 7508.27 | 3408.10 | 939.58 | 6.64 | 1.98 | -3.08 | 1250.35 | 319.34 |
|   | 4100.94 | 339.59 | 9190.56 | 5089.62 | 919.72 | 6.59 | 2.37 | -3.07 | 1874.36 | 492.92 |
| 6 | 5483.00 | 351.09 | 5834.10 | 351.10 | 1144.84 | -14.29 | 15.82 | -11.67 | 0.00 | 0.00 |
|   | 4257.00 | 356.52 | 8465.81 | 4208.81 | 636.13 | -15.60 | -17.87 | 11.55 | 82.27 | 11.77 |
|   | 5590.11 | 253.95 | 9781.67 | 4191.55 | 617.61 | -12.83 | 13.81 | -13.12 | 76.50 | 11.01 |
| 7 | 3650.50 | 205.05 | 5544.35 | 1893.85 | 438.29 | -5.21 | 6.62 | 0.10 | 14706.27 | 1424.57 |
|   | 3650.50 | 202.34 | 7662.00 | 4011.50 | 437.80 | -5.16 | 6.61 | 0.09 | 32375.88 | 5325.86 |
|   | 3650.50 | 205.13 | 9426.72 | 5776.22 | 438.15 | -5.21 | 6.62 | 0.10 | 47763.64 | 10553.81 |
| 8 | 4804.87 | 321.39 | 8264.05 | 3459.18 | 1141.92 | -10.24 | -12.57 | -0.14 | 409.39 | 81.06 |
|   | 4803.57 | 322.13 | 9868.32 | 5064.75 | 1143.05 | -10.18 | -12.50 | -0.10 | 611.95 | 121.36 |
| 9 | 4436.62 | 152.66 | 5547.61 | 1110.99 | 457.87 | -9.83 | -12.31 | 7.25 | 10.63 | 0.80 |
|   | 5395.86 | 170.98 | 7593.49 | 2197.64 | 855.62 | -14.83 | 22.08 | -2.71 | 53.50 | 3.72 |
|   | 5397.67 | 168.77 | 9639.38 | 4241.71 | 874.91 | -14.88 | 22.09 | -2.65 | 109.25 | 7.40 |
| 10 | 5216.36 | 780.66 | 6007.97 | 791.61 | 273.63 | -4.33 | -0.51 | -14.04 | 0.00 | 0.00 |
|    | 5652.58 | 330.37 | 8033.43 | 2380.85 | 701.13 | -8.34 | 4.49 | -11.69 | 1.98 | 0.55 |
|    | 6801.12 | 194.28 | 10058.89 | 3257.77 | 1054.79 | -10.89 | 10.10 | -11.24 | 5.84 | 1.63 |
| 11 | 3932.94 | 215.43 | 5681.67 | 1748.72 | 669.74 | -6.89 | 8.24 | 0.63 | 23841.97 | 4681.24 |
|    | 3919.51 | 222.80 | 7630.18  | 3710.67 | 680.47  | -6.99  | 8.35   | 0.85   | 57010.04  | 16165.55  |
|    | 3923.98 | 220.05 | 9578.69  | 5654.71 | 678.18  | -6.96  | 8.34   | 0.77   | 89492.17  | 32181.45  |
|    | 3650.50 | 68.57  | 5680.49  | 2029.99 | 964.32  | -9.98  | -12.47 | 3.26   | 100650.26 | 19921.59  |
| 12 | 6531.11 | 99.78  | 7624.88  | 1093.77 | 1041.12 | -9.12  | -11.02 | 2.27   | 50732.19  | 8718.36   |
|    | 3650.50 | 70.03  | 9569.27  | 5918.77 | 1029.02 | -9.61  | -12.19 | 3.23   | 308146.24 | 116366.73 |
|    | 3650.50 | 220.97 | 5969.30  | 2318.80 | 555.76  | -8.57  | 3.11   | 2.75   | 44922.04  | 15913.59  |
| 13 | 3650.50 | 219.50 | 7391.97  | 3741.47 | 564.95  | -8.56  | 3.43   | 2.65   | 76038.17  | 30907.23  |
|    | 3650.50 | 218.10 | 9525.99  | 5875.49 | 571.35  | -8.55  | 3.72   | 2.57   | 123008.50 | 58697.13  |
| 14 | 6583.79 | 136.24 | 7529.57  | 945.78  | 309.54  | -9.90  | -15.77 | -4.40  | 4.89      | 0.34      |
|    | 6581.72 | 138.43 | 9874.18  | 3292.46 | 312.38  | -9.87  | -15.70 | -4.16  | 19.67     | 1.34      |
|    | 3650.50 | 227.20 | 5583.43  | 1932.93 | 1359.20 | -9.44  | -11.63 | -11.80 | 265233.42 | 174397.23 |
| 15 | 3650.50 | 227.21 | 7397.08  | 3746.58 | 1359.50 | -9.44  | -11.63 | -11.80 | 527916.98 | 424692.29 |
|    | 3650.50 | 227.25 | 9210.74  | 5560.24 | 1360.55 | -9.44  | -11.63 | -11.79 | 790680.51 | 682933.87 |
|    | 4578.23 | 60.07  | 5867.83  | 1289.60 | 1360.86 | -8.78  | -4.69  | 2.62   | 39109.10  | 5733.22   |
| 16 | 4455.05 | 118.79 | 8298.70  | 3843.65 | 793.55  | -9.69  | 11.57  | -0.57  | 65299.94  | 27043.67  |
|    | 4546.08 | 97.26  | 9514.13  | 4968.05 | 1804.38 | -8.00  | -4.63  | 2.20   | 187881.19 | 44490.00  |
|    | 5346.70 | 212.60 | 5566.21  | 219.52  | 273.45  | -25.20 | 24.83  | -3.26  | 30.66     | 2.50      |
| 17 | 3867.85 | 411.59 | 8067.07  | 4199.22 | 407.68  | -14.41 | 17.11  | -4.25  | 76168.47  | 4546.10   |
|    | 3864.13 | 417.78 | 9317.50  | 5453.37 | 418.41  | -14.38 | 17.33  | -4.48  | 105870.84 | 6119.94   |
|    | 3755.42 | 424.18 | 5573.42  | 1818.00 | 872.68  | 4.70   | 1.17   | -1.58  | 755.27    | 219.38    |
| 18 | 4131.45 | 336.26 | 7592.06  | 3460.61 | 1496.89 | 3.19   | 1.99   | -1.31  | 1925.68   | 557.44    |
|    | 4129.64 | 338.16 | 9322.33  | 5192.69 | 1496.07 | 3.19   | 2.03   | -1.34  | 2965.54   | 863.27    |
|    | 3650.50 | 201.18 | 6716.97  | 3066.47 | 668.09  | -9.67  | 4.28   | -0.53  | 266564.21 | 38923.22  |
| 19 | 3650.50 | 203.42 | 8138.22  | 4487.72 | 684.28  | -9.63  | 3.62   | -0.52  | 408070.55 | 75442.05  |
|    | 3650.50 | 204.67 | 9559.48  | 5908.98 | 691.48  | -9.62  | 3.26   | -0.52  | 548972.26 | 123032.90 |
|    | 3846.60 | 230.17 | 5723.91  | 1877.31 | 753.46  | -2.96  | -2.30  | -8.95  | 37.72     | 13.03     |
| 20 | 3842.14 | 234.39 | 7755.70  | 3913.56 | 764.40  | -2.89  | -2.30  | -8.70  | 84.18     | 29.32     |
|    | 3865.47 | 105.14 | 9381.13  | 5515.66 | 648.50  | -2.65  | 3.22   | 7.78   | 90.89     | 37.14     |
|    | 3804.42 | 288.73 | 5629.83  | 1825.41 | 1578.00 | -3.06  | -7.03  | 0.10   | 24477.83  | 1051.04   |
| 21 | 3826.17 | 272.04 | 7737.97  | 3911.80 | 1536.69 | -3.08  | -7.13  | 0.15   | 54828.62  | 3768.38   |
|    | 4795.22 | 308.28 | 9319.08  | 4523.86 | 2103.02 | -3.17  | -2.02  | -0.16  | 96856.26  | 7416.74   |
|    | 3753.67 | 283.60 | 6936.07  | 3182.40 | 1963.39 | -9.36  | -9.93  | -0.07  | 52177.12  | 3830.89   |
| 22 | 3724.34 | 311.75 | 8417.98  | 4693.64 | 1936.20 | -9.25  | -9.76  | -0.02  | 76581.96  | 8626.96   |
|    | 3754.23 | 283.40 | 9899.90  | 6145.66 | 1954.74 | -9.39  | -9.99  | -0.07  | 103984.39 | 16160.81  |
| 23 | 4108.27 | 595.86 | 8293.93  | 4185.66 | 347.87  | -21.69 | -23.22 | -1.66  | 210.06    | 1.56      |
|    | 4108.03 | 597.02 | 10151.47 | 6043.44 | 346.45  | -21.67 | -23.09 | -1.63  | 311.42    | 2.34      |
|    | 4638.90 | 146.75 | 6274.21  | 1635.31 | 465.31  | -15.75 | 23.17  | -1.52  | 250.69    | 2.69      |
| 24 | 4639.23 | 146.04 | 7702.81  | 3063.57 | 474.87  | -15.83 | 23.23  | -1.54  | 513.35    | 5.44      |
|    | 6077.27 | 137.19 | 10560.01 | 4482.74 | 288.77  | -19.62 | 26.95  | -1.15  | 584.01    | 6.40      |
|    | 5216.27 | 357.41 | 6144.63  | 928.36  | 596.31  | -12.15 | -16.67 | -12.69 | 654.03    | 102.24    |
| 25 | 5950.28 | 307.64 | 7517.33  | 1567.04 | 735.93  | -11.77 | -15.09 | -12.39 | 1555.50   | 204.19    |
|    | 3774.30 | 434.60 | 9576.37  | 5802.07 | 303.65  | -12.35 | -19.21 | -11.18 | 2491.86   | 543.13    |
| 26 | 6044.00 | 152.45 | 8036.66  | 1992.67 | 359.50  | -19.42 | 26.10  | 8.57   | 53436.71  | 11541.4   |
|    | 6058.42 | 141.37 | 9901.02  | 3842.60 | 278.23  | -18.56 | 26.17  | 7.89   | 81765.09  | 20810.0   |
| 27 | 5844.80 | 240.64 | 7674.41  | 1829.61 | 588.33  | -11.27 | 10.11  | 2.14   | 897.61    | 103.12    |
|    | 7322.19 | 312.87 | 9245.58  | 1923.40 | 1034.88 | -8.85  | 8.22   | 2.51   | 1128.88   | 238.36    |
|    | 5483.00 | 166.08 | 5654.45  | 171.45  | 600.25  | -9.33  | 9.06   | 4.28   | 2.55      | 0.50      |
| 28 | 3921.63 | 366.06 | 8311.08  | 4389.44 | 696.64  | -7.07  | 5.40   | 4.15   | 10707.80  | 4942.53   |
|    | 3919.12 | 369.63 | 9639.39  | 5720.27 | 713.03  | -7.01  | 5.31   | 4.07   | 14551.53  | 6592.93   |
|    | 3915.75 | 391.15 | 5785.14  | 1869.39 | 1537.58 | -14.27 | 22.63  | -2.96  | 5522.96   | 50.88     |
| 29 | 3905.58 | 401.40 | 8002.33  | 4096.75 | 1527.61 | -14.46 | 22.84  | -2.82  | 13908.93  | 329.07    |
|    | 3907.39 | 399.72 | 9480.46  | 5573.07 | 1494.54 | -14.66 | 23.04  | -2.89  | 19345.68  | 641.45    |
|    | 4654.83 | 160.20 | 5857.48  | 1202.65 | 595.10  | -14.02 | 20.39  | 0.92   | 1123.71   | 47.55     |
| 30 | 4655.90 | 159.89 | 7998.45  | 3342.56 | 578.40  | -13.89 | 20.25  | 0.90   | 3261.83   | 146.24    |
|    | 4655.25 | 159.14 | 10139.43 | 5484.18 | 607.49  | -14.13 | 20.49  | 0.91   | 5743.99   | 251.66    |

The last column of Table 2 highlights the net change in the MOID, i.e. $\|\Delta \mathbf{r} + \delta \mathbf{r}\| - \|\Delta \mathbf{r}\|$; the table shows that the value of the deviation $\delta r$ can be significantly higher than the actual modification of the MOID. Furthermore, from the comparison between the actual achieved deviations of asteroids with small and big $\Delta r$ we can infer that $\Delta r$ itself plays an important role and cannot be neglected when dealing with a realistic impact scenario.

It should be noted that a number of solutions computed with the single objective approach belong to the set of Pareto optimal solutions. As an example, Fig. 21 shows the Pareto optimal set for asteroid 2002VU17; the black points represent the solutions within an upper boundary of 5 years, the dark gray points the solutions with an upper boundary of 10 years and the light gray points the solution with an upper boundary of 15 years. The three circles represent the Pareto optimal solutions with the maximum deviation of each scenario; the three crosses instead are the optimal solutions from the single-objective optimization with the maximum deviation of each scenario.

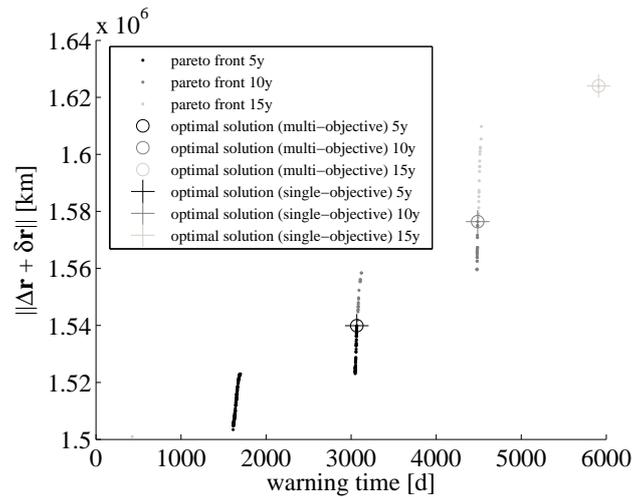

**Fig. 21    Pareto front for asteroid 2002VU17.**

Finally Fig. 22 represents the distribution, for all the asteroids and for both analyses, of the components of the impact velocity in the orbit plane.

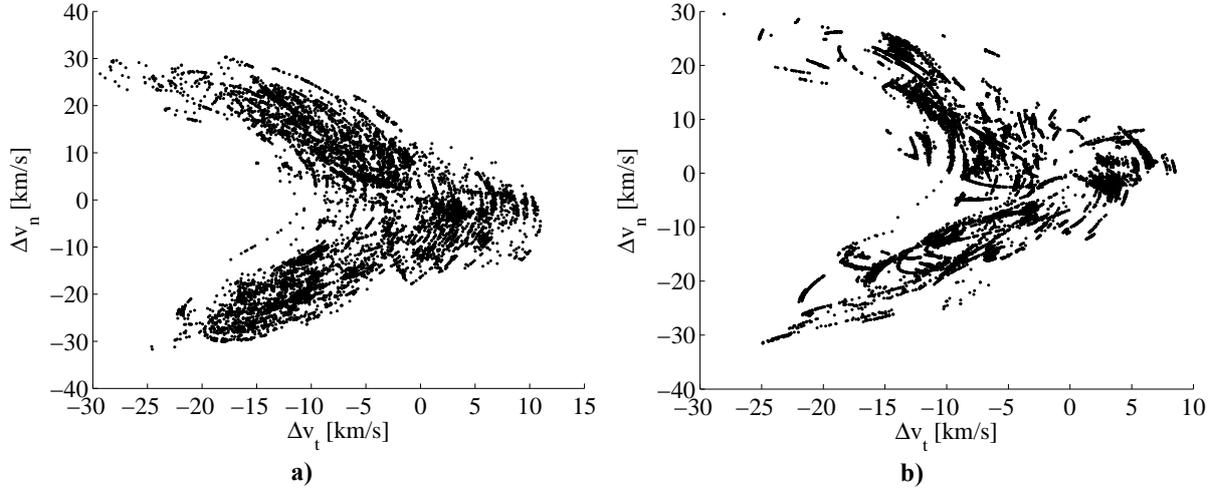

**Fig. 22** Optimal impact $\Delta v_s$ distribution for direct impacts. a) Results of the single-objective optimization. b) Results of the multi-objective optimization.

*3. Venus swing-by*

Fig. 22 shows that very high speed fuel-efficient impacts have both a very high normal and tangential components with negative sign. High speed direct impacts, therefore, correspond to trajectories that intersect almost perpendicularly the orbit of the NEO and not necessarily at the perihelion. Since this particular behavior is due to the limitation on the propellant consumption, one or more gravity assist maneuvers could improve the impact performance. Here, only the effect of a single swing-by of Venus will be considered. Table 3 reports all the solutions that show a significant improvement of the total deviation with respect to the direct transfer options.

**Table 3** Optimal launch opportunities for transfers to selected asteroids via a single Venus swing-by. Result of the single-objective optimization. The designation numbers correspond to the asteroids listed in Table 1.

| Id | $t_0$ [d] | ToF [d] | $t_{MOID}$ [MJD since 2000] | $t_w$ [d] | $m_s(t_d)$ [d] | $\Delta v_t$ [km/s] | $\Delta v_n$ [km/s] | $\Delta v_h$ [km/s] | $\delta r$ [km] | $\lVert\Delta\mathbf{r}+\delta\mathbf{r}\rVert - \lVert\Delta\mathbf{r}\rVert$ [km] |
|---|---|---|---|---|---|---|---|---|---|---|
|   | 4516.93 | 247.64  | 5618.12 | 1101.19 | 609.9125 | -4.40  | 11.85  | -3.49 | 32.26    | **9.43** |
| 4 | 3881.90 | 1292.77 | 7658.52 | 3776.62 | 551.478  | -4.86  | 11.88  | -1.28 | 99.84    | **28.77** |
|   | 3866.34 | 912.49  | 9290.84 | 5424.50 | 755.4222 | -4.69  | 11.01  | -2.60 | 252.62   | **57.18** |
|   | 3901.67 | 648.33  | 5547.61 | 1645.93 | 508.2309 | -13.31 | 27.14  | -4.44 | 15.25    | **1.51** |
| 9 | 3899.06 | 647.96  | 7593.49 | 3694.43 | 592.217  | -13.40 | 26.51  | -4.82 | 51.91    | **4.38** |
|   | 3901.41 | 648.51  | 9639.38 | 5737.97 | 510.1839 | -13.33 | 27.14  | -4.47 | 76.53    | **7.56** |
|   | 3865.16 | 272.27  | 5681.67 | 1816.51 | 776.6047 | -6.37  | 6.56   | 0.40  | 24957.00 | **4928.77** |
| 11| 3917.58 | 220.53  | 7630.18 | 3712.60 | 717.808  | -6.84  | 10.39  | -0.40 | 58123.89 | **16594.19** |
|   | 3918.67 | 218.93  | 9578.69 | 5660.02 | 719.2837 | -6.80  | 10.61  | -0.41 | 91166.93 | **32981.59** |
|   | 4455.02 | 704.40  | 5966.50 | 1511.49 | 534.0607 | -15.21 | -26.51 | -2.50 | 12.45    | **0.74** |
| 14| 4540.25 | 622.35  | 7529.57 | 2989.32 | 306.2155 | -24.34 | -32.83 | 1.04  | 32.57    | **1.59** |
|   | 3805.78 | 1363.04 | 9874.18 | 6068.40 | 429.3381 | -21.45 | -29.45 | -1.42 | 73.47    | **2.31** |
|    | 4443.39 | 447.77  | 5723.91 | **1280.52** | 272.7301 | -8.75  | -18.87| -8.33 | 20.13     | 6.16      |
| 20 | 3950.95 | 832.02  | 7349.34 | **3398.39** | 460.3267 | -6.46  | 3.25  | 8.88  | 73.59     | **30.44** |
|    | 5078.74 | 528.35  | 9381.13 | **4302.39** | 561.0051 | -5.88  | 6.81  | 8.25  | 124.05    | **69.37** |
| 26 | 3896.88 | 431.01  | 6172.31 | 2275.43     | 720.7552 | -13.63 | 16.04 | 13.56 | 71975.85  | **12271.51** |
|    | 3871.68 | 454.34  | 8036.66 | 4164.98     | 781.2508 | -16.51 | 22.58 | 9.53  | 185251.39 | **28983.47** |
|    | 3929.37 | 574.96  | 6103.23 | 2173.86     | 553.0605 | -14.63 | 13.82 | 1.81  | 1050.20   | **173.34** |
| 27 | 3921.14 | 589.63  | 7674.41 | 3753.26     | 689.3104 | -13.56 | 10.93 | 1.58  | 2496.01   | **333.90** |
|    | 5629.76 | 433.82  | 9245.58 | 3615.83     | 517.9952 | -10.32 | 6.13  | 3.27  | 1318.43   | 279.78    |
| 29 | 3930.58 | 1120.13 | 5785.14 | **1854.56** | 403.3459 | -24.39 | 31.21 | -1.28 | 1308.55   | 14.35     |
|    | 4457.15 | 1349.65 | 8002.33 | **3545.19** | 563.9278 | -18.86 | 30.26 | -1.23 | 5419.39   | **354.45** |

## IV.   Conclusion

In this paper a simple analytical expression based on proximal motion equations is derived for the computation of the deflection of potentially hazardous asteroids. An analysis of the accuracy of the proposed analytical formulation has shown its accuracy for a wide range of orbit geometries and for different deviation strategies. This formulation represents an extension of all the approaches based on a variation of the mean motion of the asteroid. Furthermore it is less computationally expensive than the approaches based on the use of the Lagrange coefficients. The proposed formulation was used at first to predict the optimal direction of the deflection $\delta \mathbf{v}$ that has to be applied to the NEO. The results presented in this paper are in agreement with already existing results obtained with different techniques. This confirms the correctness of the approach and the basic assumptions that were made.

Moreover a wide range of mission opportunities was analyzed through a hybrid global search method. Optimal launch options for direct transfers and for transfers via a single Venus gravity assist maneuver were identified for a selection of 30 asteroids with different orbital characteristics and different masses.

Though the assumed impacting spacecraft mass is quite small, it can be seen that remarkable deviations can be achieved with a reasonable time-to-impact by producing a small $\delta \mathbf{v}$ along track. On the other hand, for very short times-to-impact a more consistent $\delta \mathbf{v}$ is required, especially if the gravitational effects of the Earth are considered; in this case the direction of the optimal impulse depends on the time-to-impact and the orbital parameters of the asteroid. The results obtained in this paper show that the actual achievable change in the MOID can be significantly different from the deviation $\delta \mathbf{r}$ as a consequence of the modulus and direction of the MOID vector itself. Therefore the actual MOID cannot be neglected, in general. Furthermore it was shown that the ideal point of interception of the asteroid, when the transfer is considered, is not necessarily the pericenter of the orbit of the asteroid.

The importance of the transfer trajectory suggests that more complex sequences of gravity maneuvers and multi-impulse transfers may improve the results obtained in this paper. The design of more efficient transfer trajectories is currently under investigation and will be the subject of a future publication.

## Acknowledgments

The authors would like to thank Dr. Paolo De Pascale of the European Space Operation Centre for his help and suggestions in the selection of the asteroids.